\documentclass[10pt]{article}
\usepackage{amsfonts,mathrsfs,amssymb,amsthm,mathptm}
\usepackage{amsmath,amscd}
\usepackage{mathptm,pslatex}
\usepackage{color}
\usepackage[dvips]{graphicx}
\usepackage[all]{xy}
\usepackage{graphicx}
\usepackage{fancyhdr}

\begin{document}
\sloppy
\newcommand{\dickebox}{{\vrule height5pt width5pt depth0pt}}
\newtheorem{Def}{Definition}[section]
\newtheorem{Bsp}{Example}[section]
\newtheorem{Prop}[Def]{Proposition}
\newtheorem{Theo}[Def]{Theorem}
\newtheorem{Rem}[Def]{Remark}
\newtheorem{Lem}[Def]{Lemma}
\newtheorem{Koro}[Def]{Corollary}
\newtheorem{Claim}[Def]{Claim}

\newcommand{\cpx}[1]{#1^{\bullet}}
\newcommand{\lra}{\longrightarrow}
\newcommand{\lraf}[1]{\stackrel{#1}{\lra}}
 \newcommand{\ra}{\rightarrow}
\newcommand{\F}{\mathcal {F}}
\newcommand{\Hom}{{\rm Hom}}
\newcommand{\End}{{\rm End}}
\newcommand{\Ext}{{\rm Ext}}
\newcommand{\Tor}{{\rm Tor}}
\newcommand{\pd}{{\rm proj.dim}}
\newcommand{\inj}{{\rm inj}}
\newcommand{\lgd}{{l.{\rm gl.dim}}}
\newcommand{\gld}{{\rm gl.dim}}
\newcommand{\fd}{{\rm fin.dim}}
\newcommand{\Fd}{{\rm Fin.dim}}
\newcommand{\lfd}{l.{\rm Fin.dim}}
\newcommand{\rfd}{r.{{\rm Fin.dim}}}
\newcommand{\Mod}{{\rm Mod}}
\newcommand{\opp}{^{\rm op}}
\newcommand{\Proj}{{\rm Proj}}
\newcommand{\modcat}[1]{#1\mbox{{\rm -mod}}}
\newcommand{\pmodcat}[1]{#1\mbox{{\rm -proj}}}
\newcommand{\Pmodcat}[1]{#1\mbox{{\rm -Proj}}}
\newcommand{\injmodcat}[1]{#1\mbox{{\rm -inj}}}
\newcommand{\E}{{\rm E}_{\mathcal {F}}^{{\rm F},\Phi}}
\newcommand{\X}{ \mathscr{X}_{\mathcal {F}}^{{\rm F},\Phi}}
\newcommand{\Y}{\mathscr{Y}_{\mathcal {F}}^{{\rm F},\Phi}}
\newcommand{\A}{\mathcal {A}}
\newcommand{\C}{\mathscr{C}}
\newcommand{\D}{\mathscr{D}}
\newcommand{\h}{\mathscr{H}}
\newcommand{\K}{\mathscr{K}}
\newcommand{\Db}[1]{\mathscr{D}^{\rm b}(#1)}
\newcommand{\Cb}[1]{\mathscr{C}^b(#1)}
\newcommand{\Kb}[1]{\mathscr{K}^b(#1)}
\newcommand{\otimesL}{\otimes^{\rm\bf L}}
\newcommand{\otimesP}{\otimes^{\bullet}}
\newcommand{\rHom}{{\rm\bf R}{\rm Hom}}
\newcommand{\projdim}{\pd}
\newcommand{\stmodcat}[1]{#1\mbox{{\rm -{\underline{mod}}}}}
\newcommand{\Modcat}[1]{#1\mbox{{\rm -Mod}}}
\newcommand{\procat}[1]{#1\mbox{{\rm -proj}}}
\newcommand{\Tr}{{\rm Tr}}
\newcommand{\add}{{\rm add}}
\newcommand{\Add}{{\rm Add}}
\newcommand{\I}{{\rm Im}}
\newcommand{\Ker}{{\rm Ker}}
\newcommand{\EA}{\rm E^\Phi_\mathcal {A}}
\newcommand{\pro}{{\rm pro}}
\newcommand{\Coker}{{\rm Coker}}
\newcommand{\id}{{\rm id}}
\renewcommand{\labelenumi}{\Alph{enumi}}
\newcommand{\M}{\mathcal {M}}
\newcommand{\Mf}{\rm \mathcal {M}^f}
\newcommand{\rad}{{\rm rad}}
\newcommand{\injdim}{{\rm inj.dim}}
\newcommand{\soc}{{\rm soc}}

{\Large \bf
\begin{center}
Algebraic $K$-theory of endomorphism rings
\end{center}}
\medskip

\centerline{{\bf Hongxing Chen} and {\bf Changchang Xi$^*$}}

\renewcommand{\thefootnote}{\alph{footnote}}
\setcounter{footnote}{-1} \footnote{ $^*$ Corresponding author.
Email: xicc@bnu.edu.cn; Fax: 0086 10 58808202; Tel.: 0086 10
58808877.}
\renewcommand{\thefootnote}{\alph{footnote}}
\setcounter{footnote}{-1} 
\renewcommand{\thefootnote}{\alph{footnote}}
\setcounter{footnote}{-1} 

\begin{abstract}
We establish formulas for computation of the higher algebraic
$K$-groups of the endomorphism rings of objects linked by a morphism
in an additive category. Let ${\mathcal C}$ be an additive category,
and let $Y\ra X$ be a covariant morphism of objects in ${\mathcal
C}$. Then $K_n\big(\End_{\mathcal C}(X\oplus Y)\big)\simeq
K_n\big(\End_{{\mathcal C},Y}(X)\big)\oplus K_n\big(\End_{\mathcal
C}(Y)\big)$ for all $1\le n\in \mathbb{N}$, where $\End_{{\mathcal
C},Y}(X)$ is the quotient ring of the endomorphism ring
$\End_{\mathcal C}(X)$ of $X$ modulo the ideal generated by all
those endomorphisms of $X$ which factorize through $Y$. Moreover,
let $R$ be a ring with identity, and let $e$ be an idempotent
element in $R$. If $J:=ReR$ is homological and $_RJ$ has a finite
projective resolution by finitely generated projective $R$-modules,
then $K_n(R)\simeq K_n(R/J)\oplus K_n(eRe)$ for all $n\in
\mathbb{N}$. This reduces calculations of the higher algebraic
$K$-groups of $R$ to those of the quotient ring $R/J$ and the corner
ring $eRe$, and can be applied to a large variety of rings:
Standardly stratified rings, hereditary orders, affine cellular
algebras and extended affine Hecke algebras of type $\tilde{A}$.
\end{abstract}


\section{Introduction}
Algebraic $K$-theory collects elaborate invariants for rings. One of
the most fundamental and important questions in this theory is to
understand and calculate these invariants: algebraic $K$-groups
$K_n$ of rings. Unfortunately, this question is so hard that, up to
now,  only a few rings have gotten their algebraic $K$-groups
satisfactorily calculated (for example, see \cite{qff, rosenberg}
for details), though general, abstract algebraic $K$-theories have
been explosively developed in the last a few decades. Thus, it
becomes more reasonable to look at relationship between algebraic
$K$-groups of different rings linked by certain equivalences,
homomorphisms, or functors between their relevant categories (for
example, see \cite{DS, tt, nr}). In this way, one may compute the
algebraic $K$-groups of a ring through those of another ring.
Hopefully, this might lead to some knowledge on comprehensive
understanding of higher algebraic $K$-groups $K_n$ for rings.

This paper is a continuation of the project in this direction
started in \cite{x1} where the techniques of derived equivalences
developed in the representation theory of algebras were used to
calculate the algebraic $K$-groups of rings. More precisely, we
employed $\mathcal{D}$-split sequences introduced in \cite{hx2} to
give formulas for computation of the higher algebraic $K$-groups of
a class of rings including many maximal orders in noncommutative
algebraic geometry and in arithmetical representations (see
\cite{ch2, cr}). One of the key techniques used in \cite{x1} is to
embed a given ring into a big ring which is projective as a module
over the given ring. This method is powerful for the rings
considered there. But, for a general ring, we do not know the
existence of such an embedding. For instance, let $A$ be a
commutative ring with $I$ an ideal in $A$, we do not know how to
embed the matrix ring $R:=\begin{pmatrix}
  A             &  I    \\
  I    & A        \\
  \end{pmatrix}$
into a ring $S$ such that $_RS$ is a finitely generated projective
module and that the algebraic $K$-groups of $R$ can be computed
through those of $S$. Rings of matrix form are of importance
because, for instance, they are the essential ingredients in the
study of canonical singularities and minimal model program for
orders over surfaces (see \cite{ch2}), and of Hecke orders for
integral representations of finite groups (see \cite{Rogg, cps}).
Thus, it would be interesting to know the $K$-theory of this kind of
rings. This motivates us to consider the following general question
(see \cite{x1}):

{\bf Question.} {\it Suppose that $I$ and $J$ are two arbitrary
ideals in a ring $R$ with identity. For the ring $S:=\begin{pmatrix}
R & I  \\
J& R \end{pmatrix}$, can one give a formula for $K_n(S)$ in terms of
$K_n$-groups of quotient rings of $R$ by ideals produced from $I$
and $J$} ?

\medskip
In this paper, we will consider, more generally, algebraic
$K$-theory of the endomorphism rings of objects in a additive
category, which are linked by a morphism. As a byproduct of our
consideration, we get partial answers to the above-mentioned
question. Our idea in this paper is again to use
representation-theoretic methods for developing general results for
calculations of the higher algebraic $K$-groups of rings which
include particularly some rings mentioned above and cover also many
other interesting classes of rings such as standardly stratified
rings, hereditary orders, affine cellular algebras and extended
affine Hecke algebras of type $\tilde{A}$ (see \cite{cps, cr,
lusztig, kx}). Two key ingredients of our proofs in this paper are
recollements of triangulated categories in \cite{BBD} and the
Thomason-Waldhausen Localization Theorem which is due to Thomason
\cite[1.9.8, 1.8.2]{tt} based on the work of Waldhausen \cite{wald}.

Before stating our main results, we first introduce the notion of
covariant morphisms in an additive category (see Section \ref{sect3}
for more details).

Let $\mathcal C$ be an additive category, and let $X, Y$ be objects
in $\mathcal C$. A morphism $\lambda: Y\ra X$ in $\mathcal C$ is
said to be \emph{$X$-covariant} if the induced map $\Hom_{\mathcal
C}(X, \lambda): \Hom_{\mathcal C}(X,Y)\lra \Hom_{\mathcal C}(X,X)$
is a split monomorphism of $\End_{\mathcal C}(X)$-modules; and
\emph{covariant} if the induced map $\Hom_{\mathcal C}(X,\lambda):
\Hom_{\mathcal C}(X,Y)\lra \Hom_{\mathcal C}(X,X)$ is injective and
the induced map $\Hom_{\mathcal C}(Y,\lambda): \Hom_{\mathcal
C}(Y,Y)\lra \Hom_{\mathcal C}(Y,X)$ is a split epimorphism of
$\End_{\mathcal C}(Y)$-modules.  For example, if $\mathcal C$ is the
module category of a unitary ring $R$ and if $X$ is an $R$-module,
then, for every submodule $Y$ of $X$ with $\Hom_R(Y,X/Y)=0$, the
inclusion map is covariant, and for an idempotent ideal $I$ in $R$,
the inclusion from $I$ into $R$ is also covariant. Note that
covariant homomorphisms arise also from Auslander-Reiten sequences
and GV-ideals.

Let $\End_{{\mathcal C},Y}(X)$ denote the quotient ring of the
endomorphism ring $\End_{\mathcal C}(X)$ of the object $X$ modulo
the ideal generated by all those endomorphisms of $X$ which
factorize through the object $Y$. For example, if $I$ is an ideal in
a ring $R$ with identity, then $\End_{R,I}(R)\simeq R/I$.

Recall that an ideal $J$ in a ring $R$ is said to be
\emph{homological} if $\Tor_j^R(R/J,R/J)=0$ for all $j>0$, and that
an $R$-module $M$ has a \emph{finite-type resolution}  if it has a
finite projective resolution by finitely generated projective
$R$-modules, that is, there is an exact sequence $0\ra P_n\ra
\cdots\ra P_1\ra P_0\ra M\ra 0$ for some $n\in \mathbb{N}$ such that
all $R$-modules $P_j$ are projective and finitely generated.

Throughout this paper, we denote by $K_n(R)$ the $n$-th algebraic
$K$-group of a ring $R$ in the sense of Quillen.

Our main results in this paper are the following theorems in which
Theorem \ref{thm1} is, in some sense, a replacement of the excision
theorem of algebraic $K$-theory of rings with idempotent ideals. That is, 
it establishes a relationship between algebraic $K$-groups of rings linked by a special surjection.

\begin{Theo}  Let $R$ be a ring with identity, and $I$ an ideal of $R$.

$(1)$ If $I^2=I$, then the $K$-theory space of $\End_R(R\oplus I)$
is homotopy equivalent to the product of the $K$-theory spaces of
$\End_R(R/I)$ and $\End_{R}(I)$, and therefore
$K_n\big(\End_R(R\oplus I)\big)\simeq K_n(R/I)\oplus
K_n\big(\End_R(I)\big)$ for all $n\in \mathbb{N}$. In particular, if
the idempotent ideal $I$ is projective and finitely generated as a
left $R$-module, then  $K_n(R)\simeq K_n(R/I)\oplus
K_n\big(\End_R(I)\big)$ for all $n\in \mathbb{N}$.

$(2)$ If  $I=ReR$ for $e^2=e\in R$ such that $I$ is homological and
has a finite-type resolution as a left $R$-module, then the K-theory
space of $R$ is homotopy equivalent to the product of the $K$-theory
spaces of $eRe$ and $R/I$, and therefore $K_n(R)\simeq
K_n\big(\End_R(eRe)\big)\oplus K_n(R/I)$ for all $n\in \mathbb{N}$.
\label{thm1}
\end{Theo}

In Theorem \ref{thm1}, if we assume instead all conditions for right
$R$-modules (for example, in Theorem \ref{thm1} (1), assume that
$I_R$ is a finitely generated projective right $R$-module), then the
conclusion is still true because $K_n(R)\simeq K_n(R\opp)$ for all
$n\in \mathbb{N}$, where $R\opp$ is the opposite ring of $R$.

One cannot replace ``projective" in the second statement of Theorem
\ref{thm1} (1) with ``of finite projective dimension". Also, the
condition that $I$ is idempotent in Theorem \ref{thm1} (1) cannot be
dropped (see the examples in the last section).

As a consequence of Theorem \ref{thm1}, we have the following
result.

\begin{Theo}  Let $\mathcal C$ be an additive category and $f: Y\ra X$ a morphism of objects in
$\mathcal C$.

$(1)$ If $f$ is covariant, then the $K$-theory space of
$\End_{\mathcal C}(X\oplus Y)$ is homotopy equivalent to the product
of the $K$-theory spaces of $\End_{{\mathcal C},Y}(X)$ and
$\End_{\mathcal C}(Y)$. In particular, $K_n\big(\End_{\mathcal
C}(X\oplus Y)\big)\simeq K_n\big(\End_{{\mathcal C},Y}(X)\big)\oplus
K_n\big(\End_{\mathcal C}(Y)\big)$ for all $n\in \mathbb{N}$.

$(2)$ If $f$ is $X$-covariant, then the $K$-theory space of
$\End_{\mathcal C}(X\oplus Y)$ is homotopy equivalent to the product
of the $K$-theory spaces of $\End_{\mathcal C}(X)$ and
$\End_{{\mathcal C},X}(Y)$. In particular, $K_n\big(\End_{\mathcal
C}(X\oplus Y)\big)\simeq K_n\big(\End_{\mathcal C}(X)\big)\oplus
K_n\big(\End_{{\mathcal C},X}(Y)\big)$ for all $n\in
\mathbb{N}$.\label{mainthm}
\end{Theo}

For the dual statement of Theorem \ref{mainthm}, we refer the reader
to Theorem \ref{dualmainthm} in Section \ref{sect3}.

As is known, the excision theorem in algebraic $K$-theory of rings
gives a relationship of algebraic $K$-groups for rings linked by a
surjective ring homomorphism (see \cite{sw}). Similarly, Theorem
\ref{mainthm} describes a relationship of the algebraic $K$-groups
of the endomorphism rings of objects linked by a morphism in an
additive category.

Clearly, we can apply Theorem \ref{thm1} to standardly stratified
rings (see Section \ref{sect4} for definition) and get a reduction
formula for algebraic $K$-groups of this class of rings. It is worth
noting that ideals with the property mentioned in Theorem \ref{thm1}
(1) occur also frequently in matrix subrings.

As an application of Theorem \ref{thm1}, we have the following
corollary which provides a partial answer to the above question and
extends \cite[Theorem 1.1 (1)]{x1}.

\begin{Koro} Let $R$ be a ring with
identity, and let $J$ and $I_{ij}$ with $1\le i<j\le n$ be arbitrary
ideals of $R$ such that $I_{i j+1}J\subseteq I_{i\, j}$, $JI_{i\,
j}\subseteq I_{i+1 j}$ and $I_{ij}I_{jk}\subseteq I_{ik}$ for $j<k\le n$.
Define a ring
$$S:= {\begin{pmatrix}
  R & I_{12} &   \cdots &   \cdots & I_{1n}\\
  J & R   &   \ddots  & \ddots &   \vdots\\
  J^2     & J  & \ddots  & \ddots  &   \vdots  \\
  \vdots  & \ddots &\ddots & R &  I_{n-1\, n}\\
  J^{n-1} &  \cdots & J^2 & J  &R\\
\end{pmatrix}.}_{n\times n}$$ If $_RJ$ is projective and finitely generated, then
$$K_*(S)\simeq  K_*(R)\oplus \bigoplus_{j=1}^{n-1}K_*(R/I_{j\;j+1}J).$$ \label{cor3'}
\end{Koro}

Note that rings of the form in Corollary \ref{cor3'} not only cover
some of tiled orders, Hecke orders, and minimal model program for
orders over surfaces (see \cite{cr, Rogg, ch2}), but also occur in
commutative rings (see \cite[Section 7]{x1}) and stratification of
derived module categories arising from infinitely generated tilting
modules over tame hereditary algebras (see \cite{xc2}).

Theorem \ref{thm1} can also be applied to affine cellular algebras
(see \cite{kx}) and reduces their algebraic $K$-theory to the one of
affine commutative rings. In particular, we have the following
corollary about the algebraic $K$-groups of extended affine Hecke
algebras of type $\tilde{A}$ (for definition, see Section
\ref{sect4.2}).

\begin{Koro} Let $k$ be a field of characteristic $0$ and $q\in k$ such that $\sum_{w\in W_0}q^{\ell(w)}\ne 0$,
where $W_0$ is the symmetric group of $n$ letters and $\ell(w)$ is
the usual length function on $W_0$. For the extended affine Hecke
algebra $\h_k(n,q)$  of type $\widetilde{A}_{n-1}$, we have
$$ K_*(\h_k(n,q)) \simeq \bigoplus_{\bf c}K_*(R_{{\bf c}})$$for
$*\in \mathbb{N}$, where ${\bf c}$ runs over all two-sided cells of
the extended affine Weyl groups $W$ of type $\widetilde{A}_{n-1}$,
and $R_{\bf c}$ stands for the representation ring associated with
${\bf c}$. \label{cor4'}
\end{Koro}

Note that the ring $R_{\bf c}$ is a tensor product of rings of the
form $\mathbb{Z}[X_1, \cdots, X_s, X_{s+1},X_{s+1}^{-1}]$ with $s$ a
suitable natural number (see \cite{kx} for details). So the
$K$-theory of $R_{\bf c}$ is closely related to that of $\mathbb{Z}$
by the Fundamental Theorem in $K$-theory.

The contents of this paper are outlined as follows. In Section
\ref{sect2}, we provide necessary materials needed in our proofs.
For instance, we recall the Thomason-Waldhausen Localization Theorem
and the notion of recollements of triangulated categories. In
Section \ref{sect3}, we introduce the notion of covariant and
contravariant morphisms in an additive category, provide some of its
basic properties, and prove the main results, Theorem \ref{thm1} and
Theorem \ref{mainthm}. In Section \ref{sect4}, we   apply our
results to give formulas for calculations of algebraic $K$-groups
$K_n$ of some classes of rings, including standardly stratified
rings, matrix subrings, quantum Schur algebras, affine cellular
algebras, extended affine Hecke algebras of type $\tilde{A}$, and
skew group rings. This also proves Corollaries \ref{cor3'} and
\ref{cor4'}. In the last section, we exhibit a few examples to
demonstrate that some conditions of our results cannot be removed or
weakened. Also, two open questions are proposed at the end of this
section.

\section{Preliminaries\label{sect2}}

Given a ring $R$ with identity, we denote by $\Modcat{R}$ the
category of all left $R$-modules, and by $\modcat{R}$ the category
of all finitely generated left $R$-modules. As usual, by
$\pmodcat{A}$ we denote the category of all finitely generated
projective left $R$-modules. The complex, homotopy and derived
categories of $\Modcat{R}$ are denoted by $\C(\Modcat{R}),
\K(\Modcat{R})$ and $\D(\Modcat{R})$, respectively.

The category $\modcat{R}$ with short exact sequences is an exact
category in the sense of Quillen (see \cite{Quillen}), we denote its
$K$-theory by $G_*(R)$. As usual, we denote  by $K_*(R)$ the
$K$-theory of $\pmodcat{R}$ with split exact sequences. If $R$ is
left noetherian and has finite global dimension, then $K_*(R)\simeq
G_*(R)$ for all $*\in \mathbb{N}$. In general, even for finite
dimensional algebras over a field, the $G$-theory and $K$-theory are
not isomorphic, though the former is reduced to the one of the
endomorphism rings of simple modules.

\subsection{Waldhausen categories}

Now we recall some elementary notion about the $K$-theory of small
Waldhausen categories (see \cite{wald, tt}).

By \emph{a category  with cofibrations} we mean a category $\mathcal
C$ with a zero object $0$, together with a chosen class
co($\mathcal{C}$) of morphisms in $\mathcal{C}$ satisfying the
following three axioms:

(1) Any isomorphism in $\mathcal C$ is a morphism in
co($\mathcal{C}$),

(2) For any object $A$ in $\mathcal{C}$, the unique morphism $0\ra
A$ is in co($\mathcal{C}$),

(3) If $X\ra Y$ is a morphism in co($\mathcal{C}$), and $X\ra Z$ is
a morphism in $\mathcal{C}$, then the push-out $Y\cup_XZ$ exists in
$\mathcal{C}$, and the canonical morphism $Z\ra Y\cup_XZ$ is in
co($\mathcal{C}$). In particular. finite coproducts exist in
$\mathcal{C}$.

A morphism in co($\mathcal{C}$) is called a \emph{cofibration}.

A category $\mathcal{C}$ with cofibrations is called a
\emph{Waldhausen category} if $\mathcal{C}$ admits a class
w($\mathcal{C}$) of morphisms satisfying the two axioms:

(1) Any isomorphism in $\mathcal{C}$ is a morphism in
w($\mathcal{C}$).

(2) Given a commutative diagram
$$\begin{CD} B@<<< A@>>> C\\
@V{}VV @V{}VV @V{}VV\\
B'@<<< A'@>>>C'
\end{CD} $$
in $\mathcal{C}$ with two morphisms $A\ra B$ and $A'\ra B'$ being
cofibrations, and with $B\ra B'$, $A\ra A'$ and $C\ra C'$ being in
w($\mathcal{C}$), then the induced morphism $B\cup_AC\lra
B'\cup_{A'} C'$ is in w($\mathcal{C}$).

The morphisms in w($\mathcal{C}$) are called \emph{weak
equivalences}. Thus a Waldhausen category consists of the triple
data: A category, cofibrations and weak equivalences.

A functor between Waldhausen categories is called an \emph{exact}
functor if it preserves zero, cofibrations, weak equivalence classes
and the pushouts along the cofibrations.

A typical example of Waldhausen categories can be obtained from
complexes of modules over rings in the following manner:

Let $R$ be a ring with identity. Let $\Cb{\pmodcat{R}}$ be the small
category consisting of all bounded complexes of finitely generated
projective $R$-modules. This is a Waldhausen category. That is, the
weak equivalences are the homotopy equivalences, and the
cofibrations are the degreewise split monomorphisms. By just
inverting the weak equivalences, we get the derived category of
$\Cb{\pmodcat{R}}$, which is the homotopy category
$\Kb{\pmodcat{R}}$ of $\Cb{\pmodcat{R}}$.

For a small Waldhausen category $\mathcal{C}$, a $K$-theory
$K_n(\mathcal{C})$ was defined in \cite{wald}. In particular, for
the small, Waldhausen category $\Cb{\pmodcat{R}}$, it is shown by a
theorem of Gillet-Waldhausen that its $K$-theory is the same as the
$K$-theory of $R$ defined by Quillen. That is, $K_n(R)\simeq
K_n(\Cb{\pmodcat{R}})$ for all $n\ge 0$.

In this paper, we always assume that all Waldhausen categories
considered are of this type, that is, they are full subcategories of
the category of chain complexes over some abelian category.

The following result, which is called the Thomason-Waldhausen
Localization Theorem in the literature, says that we can get an
exact sequence of $K$-groups of Waldhausen categories from a short
exact sequence of their derived categories, which is induced from
exact functors between the given Waldhausen categories (for example,
see \cite{tt}, \cite[Theorem 2.3]{nr}).

\begin{Theo} (Thomason-Waldhausen Localization Theorem) Let $\mathcal{R, S}$ and
$\mathcal{T}$ be small, Waldhausen categories. Suppose
$\mathcal{R}\ra \mathcal{S}\ra\mathcal{T}$ are exact functors of
Waldhausen categories. Suppose further that

$(i)$ the induced triangulated functors of derived categories
$\D(\mathcal{R})\lra \D(\mathcal{S})\lra\D(\mathcal{T})$ compose to
zero.

$(ii)$ The functor $\varphi: \D(\mathcal{R})\lra \D(\mathcal{S})$ is
fully faithful.

$(iii)$ If $x$ and $x'$ are objects of $\D(\mathcal{S})$, and the
direct sum $x\oplus x'$ is isomorphic in $\D(\mathcal{S})$ to
$\varphi(z)$ for some $z \in \D(\mathcal{R})$, then $x, x'$ are
isomorphic to $\varphi(y), \varphi(y')$ for some $y, y'\in
\D(\mathcal{R})$.

$(iv)$ The natural map $\D(\mathcal{S})/\D(\mathcal{R})\lra
\D(\mathcal{T})$ is an equivalence of categories.

Then the sequence of spectra $K(\mathcal{R})\ra K(\mathcal{S}) \ra
K(\mathcal{T})$ is a homotopy fibration, and therefore there is a
long exact sequence of $K$-theory
$$\cdots \lra K_{n+1}(\mathcal{T})\lra K_n(\mathcal{R})\lra K_n(\mathcal{S})\lra
K_n(\mathcal{T}) \lra K_{n-1}(\mathcal{R})\lra$$ $$\cdots\lra
K_0(\mathcal{S})\lra K_0(\mathcal{R})\lra K_0(\mathcal{T})$$ for all
$n\in \mathbb{N}.$ \label{waldhausen}
\end{Theo}

\subsection{Recollements}
Another notion needed in our proofs is recollements which were
introduced by Beilinson, Bernstein and Deligne (see \cite{BBD}) to
study the behaves of triangulated category of perverse sheaves of
geometric objects.

Let $\cal D$ be a triangulated categories with a shift functor
denoted by [1]. 

Let $\mathcal{D'}$ and $\mathcal{D''}$ be triangulated categories.
We say that $\mathcal{D}$ is a \emph{recollement} of $\mathcal{D'}$
and $\mathcal{D''}$ if there are six triangle functors as in the
following diagram
$$\xymatrix{\mathcal{D''}\ar^-{i_*=i_!}[r]&\mathcal{D}\ar^-{j^!=j^*}[r]
\ar^-{i^!}@/^1.2pc/[l]\ar_-{i^*}@/_1.6pc/[l]
&\mathcal{D'}\ar^-{j_*}@/^1.2pc/[l]\ar_-{j_!}@/_1.6pc/[l]}$$ such
that

$(1)$ $(i^*,i_*),(i_!,i^!),(j_!,j^!)$ and $(j^*,j_*)$ are adjoint
pairs;

$(2)$ $i_*,j_*$ and $j_!$ are fully faithful functors;

$(3)$ $i^!j_*=0$ (and thus also $j^! i_!=0$ and $i^*j_!=0$); and

$(4)$ for each object $X\in\mathcal{D}$, there are two canonical
triangles in $\mathcal D$:
$$
i_!i^!(X)\lra X\lra j_*j^*(X)\lra i_!i^!(X)[1],
$$
$$
j_!j^!(X)\lra X\lra i_*i^*(X)\lra j_!j^!(X)[1],$$ where
$i_!i^!(X)\ra X$ and $j_!j^!(X)\ra X$ are counit adjunction maps,
and where $X\ra j_*j^*(X)$ and $ X\ra i_*i^*(X)$ are unit adjunction
maps.

It is know from the definition of recollements that the Verdier
quotients of $\mathcal{D}$ by the images of the triangle functors
$j_!$ and $i_*$ are equivalent to $\mathcal{D''}$ and
$\mathcal{D'}$, respectively.

A typical example of recollements occurs in the following situation.
Let $R$ be a ring with an idempotent ideal $J=ReR$ for $e^2=e\in R$.
Suppose that $J$ is a stratifying ideal of $R$ (for definition, see
Section \ref{sect4}), then there is a recollement:

$$\xymatrix{\D(\Modcat{R/J})\ar^-{D(\lambda_*)}[r]&\D(\Modcat{R})\ar^-{eR\otimes_R^{\mathbb{L}}-}[r]
\ar^-{\mathbb{R}\Hom_R(R/J,-)\;}@/^1.2pc/[l]\ar_-{R/J\otimes_R^{\mathbb
L}-}@/_1.6pc/[l] &\D(\Modcat{eRe})\ar_-{Re\otimes_{eRe}^{\mathbb
L}-}@/_1.6pc/[l]\ar^-{\;\mathbb{R}\Hom_{eRe}(eR,-)}@/^1.2pc/[l] }$$
where $\lambda_*: \Modcat{R/J}\lra \Modcat{R}$ is the restriction
functor, $Re\otimes_{eRe}^{\mathbb{L}}-$ is the total left-derived
functor of $Re\otimes_{eRe}-$, and $\mathbb{R}\Hom_{eRe}(eR,-)$ is
the total right-derived functor of $\Hom_{eRe}(eR,-)$. For more
details, we refer the reader to \cite{cps}.

Note that $\D(\Modcat{R})$ is a triangulated category with small
coproducts.

\medskip
Finally, we point out the following homological fact which is needed
in our proofs.

\begin{Lem} Let $R$ be a ring with identity, and let $J=ReR$ for $e^2=e\in R$. Suppose that $M$ is an
$R$-module with the following two properties:

$(1)$ $\Tor_j^R(R/J,M)=0$ for all $j\ge 0$, and

$(2)$ $M$ has a finite-type resolution, that is, there is an exact
sequence $0\lra P'_n\lraf{\epsilon'_n} \cdots \lra
P'_1\lraf{\epsilon'_1} P'_0\lraf{\epsilon'_0} M\ra 0$ with all
$P'_j$ finitely generated projective $R$-modules.

Then there is an exact sequence of $R$-modules:
$$0\ra P_n\lra \cdots \lra P_1\lra
P_0\lra M\lra 0  $$ such that all $P_j$ lie in $\add(Re)$.
\label{tor-seq}
\end{Lem}

{\it Proof.} This result is known for modules over Artin algebras,
where one may use minimal projective resolutions (see \cite{apt}).
For general rings, projective covers of modules may not exist. For
convenience of the reader, we include here a proof.

Given such a sequence in (2), we define $K_i'=\Ker(\epsilon'_{i-1})$
for $1\le i\le n$. Then $K'_i$ is finitely generated.

It follows from $\Tor_0^R(R/J,M)=0$ that $JM=M$. Since the trace of
$Re$ in the module $M$ is just $JM$ and since $M$ is finitely
generated, there is a finite index set $I_0$ and a surjective
homomorphism $ P_0:=\bigoplus_{i\in I_0}Re\lraf{\epsilon_0} M$. We
define $K_1=\Ker(\epsilon_0)$. Then, by Schanuel's Lemma, we have
$K_1\oplus P'_0\simeq K'_1\oplus P_0$, and therefore $K_1$ is
finitely generated. It follows from $\Tor_1^R(R/J,M)=0$ that the
sequence
$$ 0\lra K_1/JK_1\lra P_0/JP_0\lra M/JM\lra 0 $$
is exact. This means that $JK_1=K_1$ because $JP_0=P_0$. Observe
that $\Tor_1^R(R/J,K_1)=\Tor_2^R(R/J,M)=0$. So, for $K_1$, we can do
the similar procedure as we did above and get a surjective
homomorphism $P_1:=\bigoplus_{j\in I_1}Re\lraf{\epsilon_1} K_1$ with
$I_1$ a finite set, such that $K_2:=\Ker(\epsilon_1)$ is finitely
generated with $JK_2=K_2$ and $\Tor_1^R(R/J,K_2)=0$. Hence, by using
the generalized Schanuel's Lemma, we can iterate this procedure.
Since the projective dimension of $M$ is finite,  we must stop after
$n$ steps and reach at a desired sequence mentioned in the lemma.
$\square$

\medskip
{\it Remark.} The above proof shows that for an $R$-module $M$, the
condition (1) is equivalent to

(2') There is a projective resolution $\cdots \ra P_n\ra\cdots\ra
P_1\ra P_0\ra M\ra 0$ such that $P_j\in \Add(Re)$, where $\Add(M)$
is the full subcategory of $\Modcat{R}$ consisting of all those
$R$-modules which are direct summands of direct sums of copies of
$M$.

Thus $J=ReR$ is homological if and only if such a sequence (2') for
$_RJ$ exists.

\section{Proofs of the main results \label{sect3}}

In this section, we introduce the notion of covariant and
contravariant morphisms, and prove the main results, Theorems
\ref{thm1} and \ref{mainthm}.

Let $R$ be a ring with identity and $X$ an $R$-module. A submodule
$Y$ of $X$ is a \emph{trace} in $X$ if $\Hom_R(Y,X/Y)=0$, and a
\emph{weak trace} in $X$ if the inclusion from $Y$ to $X$ induces an
isomorphism $\Hom_R(Y,Y)\ra\Hom_R(Y,X)$ of abelian groups. For
example, every idempotent ideal of $R$ is a trace of the regular
$R$-module $_RR$, and every GV-ideal $J$ of $R$ is a weak trace of
$_RR$ (see \cite[Section 7]{x1} for definition). Also, the socle of
any finite dimensional algebra $A$ over a field is a weak trace of
$_AA$. In particular, the socle of the ring $R:=\mathbb{Q}[X]/(X^2)$
is a weak trace in $R$, but not a trace in $R$.

In general, for any $R$-modules $X$ and $Y$, there is a recipe for
getting weak trace submodules of $X$. Let $t_Y(X)$ be the sum of all
images of homomorphisms from $Y$ to $X$ of $R$-modules. Then
$t_Y(X)$ is a weak trace of $X$.

We denote by $\End_{R,Y}(X)$ the quotient ring of the endomorphism
ring $\End_R(X)$ of $X$ modulo the ideal generated by those
endomorphisms of $X$ which factorize through the module $Y$. Note
that this ideal consists actually of all those endomorphisms $f:
X\ra X$ such that there is an $R$-module $Z\in \add(Y)$ and
homomorphisms $g: X\ra Z$ and $h: Z\ra X$ with $f=gh$, where
$\add(Y)$ is the full subcategory of $R$-Mod consisting of all
modules which are direct summands of direct sums of finitely many
copies of $Y$. For instance, if $I$ is an ideal in $R$, then
$\End_{R,I}(R)\simeq R/I$.

Motivated by weak trace submodules, we introduce the following
notion.

\begin{Def} Let $\mathcal C$ be an additive category. A morphism $\lambda: Y\ra X$ of objects in $\mathcal C$ is said to be covariant
if

$(1)$ the induced map $\Hom_{\mathcal C}(X,\lambda): \Hom_{\mathcal
C}(X,Y)\ra \Hom_{\mathcal C}(X,X)$ is injective, and

$(2)$ the induced map $\Hom_{\mathcal C}(Y,\lambda): \Hom_{\mathcal
C}(Y,Y)\ra \Hom_{\mathcal C}(Y,X)$ is a split epimorphism of
$\End_{\mathcal C}(Y)$-modules.

Dually, a morphism $\beta: N \ra M$ in $\mathcal C$ is said to be
contravariant if

$(1')$ $\Hom_{\mathcal C}(\beta, N):  \Hom_{\mathcal C}(M,N)\ra
\Hom_{\mathcal C}(N,N)$ is injective, and

$(2')$ $\Hom_{\mathcal C}(\beta, M):\Hom_{\mathcal C}(M,M)\ra
\Hom_{\mathcal C}(N,M)$ is a split epimorphism of right
$\End_{\mathcal C}(M)$-modules
\end{Def}

Clearly, if $Y$ is a weak trace submodule of an $R$-module $X$, then
the inclusion map is a covariant homomorphism. Another example of
covariant homomorphisms is the following: If $0\ra Z\ra Y\lraf{g}
X\ra 0$ is an Auslander-Reiten sequence in $R$-mod with
$\Hom_R(Y,Z)=0$, then the homomorphism $g$ is covariant.

For covariant morphisms, we have the following properties.

\begin{Lem} Let $\mathcal C$ be an additive category, and let $\lambda: Y\ra X$ be a covariant morphism of
objects in $\mathcal C$. We define $\Lambda:=\End_{\mathcal
C}(X\oplus Y)$, and let $e_Y$ be the idempotent element of $\Lambda$
corresponding to the projection onto $Y$. Then

$(1)$ $_{\Lambda}\Lambda e_Y\Lambda$ is a finitely generated
projective $\Lambda$-module.

$(2)$ The composition map $\mu: \Hom_{\mathcal
C}(X,Y)\otimes_{\End_{\mathcal C}(Y)}\Hom_{\mathcal C}(Y,X)\ra
\End_{\mathcal C}(X)$ is injective. Thus the cokernel of $\mu$ is
isomorphic to $\End_{{\mathcal C},Y}(X)$. \label{lem01}
\end{Lem}

 To prove this lemma, we use the following observation.

\begin{Lem}Let $S$ be a ring with identity, and let $e$ be an idempotent
element in $S$. Then $_SSeS$ (respectively, $SeS_S$) is  projective
and finitely generated if and only if $eS(1-e)$ (respectively,
$(1-e)Se$) is  projective and finitely generated as an $eSe$-module
(respectively, a right $eSe$-module), and the multiplication map
$${\mu}:(1-e)Se\otimes_{eSe}eS(1-e)\lra (1-e)S(1-e)$$ is injective.
 \label{projective}
\end{Lem}

{\it Proof.} Suppose that $eS(1-e)$ is a finitely generated
projective $eSe$-module and that the multiplication map
$(1-e)Se\otimes_{eSe}eS(1-e)\lraf{\mu} (1-e)S(1-e)$ is injective.
Then it is easy to see that the multiplication map
$Se\otimes_{eSe}eS\lra SeS$ is an isomorphism of $S$-$S$-bimodules.
Since $eS=eSe\oplus eS(1-e)$, we know that $_SSeS$ is  projective
and finitely generated.

Conversely, suppose that $_SSeS$ is projective and finitely
generated. One the one hand, since $_SSeS$ is projective,  we can
show that the multiplication map $\mu: Se\otimes_{eSe}eS\ra SeS$ is
injective (see \cite[Statement 7]{dr}). This implies that the map
$\mu: (1-e)Se\otimes_{eSe}eS(1-e)\lra (1-e)S(1-e)$ is injective. On
the other hand, since $_SSeS$ is finitely generated, there is a
finite subset $\{x_i\mid i\in I\}$ of $S$ such that the map
$\bigoplus_{i\in I}Se\lra SeS$, defined by $(a_i)_{i\in I}\mapsto
\sum_{i\in I}a_ix_i$, is surjective. This shows that $_SSeS$ is a
direct summand of a direct sum of finitely many copies of $Se$. Thus
$eS$ is a direct summand of a free $eSe$-module of finite rank. This
implies that the $eSe$-module $eS(1-e)$ is projective and finitely
generated.

The same arguments applies to the right module $SeS_S$. $\square$

\medskip
{\it Proof of Lemma \ref{lem01}.} Clearly, $\Lambda =
\begin{pmatrix}
  \End_{\mathcal C}(X)   &  \Hom_{\mathcal C}(X,Y)    \\
 \Hom_{\mathcal C}(Y,X)    & \End_{\mathcal C}(Y)         \\
\end{pmatrix}$. Let $e:=\begin{pmatrix}
  0   &  0    \\
 0    & 1         \\
\end{pmatrix}$ and $f:=1-e$. Thus $e=e_Y$, $e\Lambda e\simeq \End_{\mathcal C}(Y)$, $f\Lambda f\simeq \End_{\mathcal C}(X)$,
$f\Lambda e\simeq \Hom_{\mathcal C}(X,Y)$ and $e\Lambda f\simeq
\Hom_{\mathcal C}(Y,X)$, where the left $\End_{\mathcal
C}(Y)$-module structure of $\Hom_{\mathcal C}(Y,X)$ is induced from
the right $\End_{\mathcal C}(Y)$-module structure of $Y$. In the
following we will often use these identifications, as was done in
\cite{xi2}. Since $\lambda$ is a covariant homomorphism, the induced
map
$$\lambda^*=\Hom_{\mathcal C}(Y,\lambda): \Hom_{\mathcal C}(Y,Y)\lra\Hom_{\mathcal C}(Y,X)$$ is a
split epimorphism of $\End_{\mathcal C}(Y)$-modules. Thus there is a
homomorphism $\gamma: \Hom_{\mathcal C}(Y,X)\ra\Hom_{\mathcal
C}(Y,Y)$ such that $\gamma\lambda^*=id$. This means that
$\Hom_{\mathcal C}(Y, X)$ is a direct summand of the regular
$\End_{\mathcal C}(Y)$-module. Thus $e\Lambda f$ is projective and
finitely generated as a left $e\Lambda e$-module since a direct
summand of a finitely generated module is finitely generated.

Now we show that the multiplication map $f\Lambda e\otimes_{e\Lambda
e}e\Lambda f\lra f\Lambda f$ is injective. This is equivalent to
showing that the composition map
$$ \mu:\quad \Hom_{\mathcal C}(X,Y)\otimes_{\End_{\mathcal C}(Y)} \Hom_{\mathcal C}(Y, X)\lra \End_{\mathcal C}(X),$$ given
by $x\otimes f\mapsto xf$ for $x\in \Hom_{\mathcal C}(X,Y)$ and
$f\in \Hom_{\mathcal C}(Y,X),$ is injective. However, the
injectivity of $\mu$ follows from the injectivity of $\Hom_{\mathcal
C}(X,\lambda): \Hom_{\mathcal C}(X,Y)\lra \Hom_{\mathcal C}(X,X)$
together with the following commutative diagram
$$\begin{CD}
\Hom_{\mathcal C}(X,Y)\otimes_{\End_{\mathcal C}(Y)}\Hom_{\mathcal C}(Y,Y)@>{\simeq}>{\mu'}> \Hom_{\mathcal C}(X,Y)\\
@VV{\Hom_{\mathcal C}(X,Y)\otimes \lambda^*}V @VV{\Hom_{\mathcal C}(X,\lambda)}V\\
\Hom_{\mathcal C}(X,Y)\otimes_{\End_{\mathcal C}(Y)}\Hom_{\mathcal
C}(Y,X)@>{\mu}>> \Hom_{\mathcal C}(X,X)
\end{CD} $$
since the bottom $\mu$ is a composite of three injective maps, that
is, $\mu=\big(\Hom_{\mathcal C}(X,Y)\otimes \gamma\big)\mu'
\,\Hom_{\mathcal C}(X,\lambda)$. Here, we use the identity
$\gamma\lambda^*=id$. Thus, by Lemma \ref{projective}, we see that
$_{\Lambda}\Lambda e\Lambda$ is a finitely generated projective
$\Lambda$-module.

Now the second statement of Lemma \ref{lem01} is also clear.
$\square$

\smallskip
Dually, for contravariant morphisms, we have the following
statement.

\begin{Lem} Let $\mathcal C$ be an additive category, and let $\lambda: Y\ra X$ be a contravariant morphism of
objects in $\mathcal C$. We define $\Lambda:=\End_{\mathcal
C}(X\oplus Y)$, and let $e_X$ be the idempotent element of $\Lambda$
corresponding to the projection onto $X$. Then

$(1)$ $\Lambda e_X\Lambda_{\Lambda}$ is a finitely generated
projective right $\Lambda$-module.

$(2)$ The composition map $\mu: \Hom_{\mathcal
C}(Y,X)\otimes_{\End_{\mathcal C}(X)}\Hom_{\mathcal C}(X,Y)\lra
\End_{\mathcal C}(Y)$ is injective. Thus the cokernel of $\mu$ is
isomorphic to $\End_{{\mathcal C},X}(Y)$. \label{lem01'}
\end{Lem}

For convenience, we introduce the following definition of
$X$-covariant morphisms. Observe that the condition in this
definition strengthens only the first and does not involve the
second condition in the definition of covariant or contravariant
morphisms.

\begin{Def}
A morphism $f: Y\ra X$ in an additive category $\mathcal C$ is said
to be

$(1)$ $X$-covariant if the induced map $Hom_{\mathcal C}(X,f)$ is a
split monomorphism of $\End_{\mathcal C}(X)$-modules.

$(2)$ $Y$-contravariant if  the induced map $Hom_{\mathcal C}(f, Y)$
is a split monomorphism of right $\End_{\mathcal C}(Y)$-modules.
\end{Def}

Clearly, if $f: Y\ra X$ is covariant, then the inclusion from
$\Ker(f)$ into $Y$ is $Y$-covariant. Dually, if $f: Y\ra X$ is
contravariant, then the canonical surjection from $X$ to $\Coker(f)$
is $X$-contravariant. For $X$-covariant and $Y$-contravariant
morphisms, we have the following properties.

\begin{Lem} Let $\mathcal C$ be an additive category, and let $\lambda: Y\ra X$ be a morphism of
objects in $\mathcal C$. We define $\Lambda:=\End_{\mathcal
C}(X\oplus Y)$, and let $e_X$ and $e_Y$ be the idempotent elements
of $\Lambda$ corresponding to the projection onto $X$ and $Y$,
respectively.

$(1)$ If $\lambda$ is $X$-covariant, then $_{\Lambda}\Lambda
e_X\Lambda$ is a finitely generated projective $\Lambda$-module. In
this case, $\Lambda/\Lambda e_X\Lambda\simeq \End_{{\mathcal
C},X}(Y)$.

$(2)$ If $\lambda$ is $Y$-contravariant, then $\Lambda
e_Y\Lambda_{\Lambda}$ is a finitely generated projective
$\Lambda$-module. In this case, $\Lambda/\Lambda e_Y\Lambda\simeq
\End_{{\mathcal C},Y}(X)$.\label{3.4"}
\end{Lem}

{\it Proof.} (1) The proof is similar to that of Lemma \ref{lem01}.
Here, we only outline its main points.

Since $\lambda^*=\Hom_{\mathcal C}(X,\lambda): \Hom_{\mathcal
C}(X,Y)\lra \Hom_{\mathcal C}(X,X)$ is a split monomorphism of
$\End_{\mathcal C}(X)$-modules, we see that

(a) $e_X\Lambda e_Y$ is a finitely generated projective $e_X\Lambda
e_X$-module, and

(b) $\Hom_{\mathcal C}(Y,X)\otimes \lambda^*: \Hom_{\mathcal
C}(Y,X)\otimes_{\End_{\mathcal C}(X)}\Hom_{\mathcal C}(X,Y)\lra
\Hom_{\mathcal C}(Y,X)\otimes_{\End_{\mathcal C}(X)}\Hom_{\mathcal
C}(X,X)$ is a split monomorphism.

To see that the multiplication map $\mu: e_Y\Lambda
e_X\otimes_{e_X\Lambda e_X}e_X\Lambda e_Y\lra e_Y\Lambda e_Y$ is
injective, we consider the following commutative diagram:

$$\begin{CD}
\Hom_{\mathcal C}(Y,X)\otimes_{\End_{\mathcal C}(X)}\Hom_{\mathcal C}(X,Y)@>>> \Hom_{\mathcal C}(Y,Y)\\
@VV{\Hom_{\mathcal C}(Y,X)\otimes \lambda^*}V @VV{\Hom_{\mathcal C}(Y,\lambda)}V\\
\Hom_{\mathcal C}(Y,X)\otimes_{\End_{\mathcal C}(X)}\Hom_{\mathcal
C}(X,X)@>{\simeq}>>\Hom_{\mathcal C}(Y,X)
\end{CD} $$
where the horizontal maps are composition maps. This implies that
$\mu$ is injective, and therefore $\Lambda/\Lambda e_X\Lambda\simeq
\End_{{\mathcal C},X}(Y)$. Now (1) follows immediately from Lemma
\ref{projective}.

(2) The proof is left to the reader. $\square$

\medskip
{\bf Proof of Theorem \ref{thm1}.}

 (1) We first show that (1) follows from (2) and Lemma \ref{lem01}.

Assume that $I$ is an idempotent ideal of $R$. Then
$\Hom_R(I,R/I)=0$ and $I$ is a trace of $_RR$. Note that
$\End_{R,I}(R)\simeq R/I$. Thus, by Lemma \ref{lem01}, the first
statement of (1) follows from (2) immediately.

Now assume further that $_RI$ is projective and finitely generated.
Then the $R$-module $R\oplus I$ is a progenerator for $\Modcat{R}$,
and therefore $R$ and $\Lambda:=\End_R(_RR\oplus I)$ are Morita
equivalent. Hence, by the first statement of (1), $K_n(R)\simeq
K_n(\Lambda)\simeq K_n\big(\End_R(I)\big)\oplus K_n(R/I)$ for all
$n\in\mathbb{N}$. This finishes the proof of Theorem \ref{thm1} (1).

(2) Now we prove (2). This is precisely the following proposition.

\begin{Prop} Let $R$ be a ring with identity, and let $e^2=e\in R$
such that $J:=ReR$ is homological and $_RJ$ has a finite-type
resolution. Then the $K$-theory space of $R$ is homotopy equivalent
to the product of the $K$-theory spaces of $eRe$ and $R/J$, and
therefore
$$K_n(R)\simeq K_n\big(\End_R(eRe)\big)\oplus K_n(R/J)$$ for all $n\in
\mathbb{N}$. \label{idempotent}
\end{Prop}

{\it Proof.} Recall that $\Cb{\pmodcat{R}}$ is the category
consisting of all bounded complexes of finitely generated projective
$R$-modules. This is a Waldhausen category. That is, the weak
equivalences are the homotopy equivalences, and the cofibrations are
the degreewise split monomorphisms. By just inverting the weak
equivalences, we get the derived category of $\Cb{\pmodcat{R}}$,
which is $\Kb{\pmodcat{R}}$.

Let $\D^c(R)$ be the full subcategory of $\D(\Modcat{R})$ consisting
of all compact objects in $\D(\Modcat{R})$. Then $\D^c(R)$ is a
triangulated subcategory. Recall that a complex $\cpx{X}\in
\D(\Modcat{R})$ is said to be \emph{compact} if $\Hom_{\D(
R\footnotesize{\rm -Mod})}(\cpx{X},-)$ commutes with coproducts in
$\D(\Modcat{R})$. It is shown in \cite[Corollary 4.4]{nr} that
$\D^c(R)$ consists of objects which are isomorphic in
$\D(\Modcat{R})$ to bounded chain complexes of finitely generated,
projective $R$-modules. Thus, any finite direct sum of compact
objects is compact, any direct summand of a compact object is
compact, and $\Kb{\pmodcat{R}}$ is equivalent to $\D^c(R)$ as
triangulated categories.

Note that if $J:=ReR$ is homological and $_RJ$ admits a finite-type
resolution, that is, there is an exact sequence
$$ 0\lra P_n\lra \cdots\lra P_0\lra {}_RReR\lra 0$$ with all $P_j$ finitely generated projective $R$-modules, then it follows from
Lemma \ref{tor-seq} that we may assume  $P_j\in \add(Re)$ for all
$j$. Thus $eP_j\in \add(_{eRe}eRe)$ and the $eRe$-module $eR$ has a
finite-type resolution. Therefore $eR$ is a compact object in
$\D(\Modcat{eRe})$. Clearly, $R/J\in \D^c(R/J)$. Now, we consider
the small Waldhausen categories $\Cb{\pmodcat{R}}$,
$\Cb{\pmodcat{eRe}}$ and $\Cb{\pmodcat{R/J}}$, and the functors
$$ Re\otimes_{eRe}-:
\Cb{\pmodcat{eRe}}\longrightarrow \Cb{\pmodcat{R}}, \quad
(R/J)\otimes_R-: \Cb{\pmodcat{R}}\longrightarrow \Cb{\pmodcat{R/J}}.
$$
Since the exact structure of these categories is the degreewise
split short exact sequences, we see that the two functors are exact.
Moreover, it is well known (for example, see \cite{cps}) that we
have a recollement

$$\xymatrix{\D(\Modcat{R/J})\ar^-{D(\lambda_*)}[r]&\D(\Modcat{R})\ar^-{eR\otimes_R^{\mathbb{L}}-}[r]
\ar^-{}@/^1.2pc/[l]\ar_-{R/J\otimes_R^{\mathbb L}-}@/_1.6pc/[l]
&\D(\Modcat{eRe})\ar_-{Re\otimes_{eRe}^{\mathbb
L}-}@/_1.6pc/[l]\ar_-{}@/^1.2pc/[l] \ar_-{}@/_1.6pc/[l]},$$ where
$\lambda_*: \Modcat{R/J}\lra \Modcat{R}$ is the restriction functor.
Note that $R/J$ and $eR$ are compact objects in $\D(\Modcat{R})$ and
$\D(\Modcat{eRe})$, respectively. Thus the exact functors
$D(\lambda_*)$ and $eR\otimesL_{R}-$ preserve compact objects. It is
known that, for a recollement, the functors
$R/J\otimes_R^{\mathbb{L}}-$ and $Re\otimes_{eRe}^{\mathbb{L}}-$
always preserves compact objects. Thus, from the above recollement
we can get the following ``half recollement" for the subcategories
of compacts objects:

$$\xymatrix{\D^c(R/J)\ar^-{D(\lambda_*)}[r]&\D^c(R)\ar^-{eR\otimes_R^{\mathbb{L}}-}[r]
\ar_-{R/J\otimes_R^{\mathbb L}-}@/_1.6pc/[l]
&\D^c(eRe)\ar_-{Re\otimes_{eRe}^{\mathbb L}-}@/_1.6pc/[l]
\ar_-{}@/_1.6pc/[l]}\vspace{0.3cm}$$ Note that $\D^c(R)$ may not
have small coproducts in general. This half recollement implies the
following commutative diagram of triangle functors:

$$\xymatrix@C=1.2cm{\D^c(R/J)
&\D^c(R)\ar[l]_-{R/J\otimes_{R}^{\mathbb L}-} & \D^c(eRe)\ar[l]_-{Re\otimes_{eRe}^{\mathbb L}-} \\
\Kb{\pmodcat{R/J}}\ar[u]_{\simeq}
&\Kb{\pmodcat{R}}\ar[l]_-{R/J\otimes_R-}\ar[u]_{\simeq}&
\Kb{\pmodcat{eRe}}\ar[l]_-{Re\otimes_{eRe}-}\ar[u]_{\simeq}
}\vspace{0.3cm}$$ This shows that the two functors in the top row
are induced from the exact functors in the bottom row. Moreover, we
have the following properties:

(1) Clearly, it follows from the half-recollement that the
composition of the two functors $Re\otimes_{eRe}^{\mathbb{L}}-$ and
$R/J\otimes_R^{\mathbb{L}}-$ is zero, that the functor
$Re\otimes_{eRe}^{\mathbb{L}}-$ is fully faithful, and that the
natural map $\D^c(R)/\D^c(eRe)\lra \D^c(R/J)$ is an equivalence of
categories. The latter follows actually from a general known fact:
If $F: \mathcal{C}\ra \mathcal{D}$ is a triangle functor which
admits a fully faithful right adjoint functor $G:\mathcal{D}\ra
\mathcal{C}$, then $F$ induces uniquely a triangle equivalence
between $\mathcal{C}/\mbox{Ker}(F)$ and $\mathcal{D}$, where
Ker$(F)$ stands for the full subcategory of $\mathcal{C}$ consisting
of all those objects $x$ such that $F(x)=0$.

(2) If $x$ and $x'$ are two objects of $\D^c(R)$, and the direct sum
$x\oplus x'$ is isomorphic in $\D^c(R)$ to
$Re\otimes_{eRe}^{\mathbb{L}}z$ for some $z\in \D^c(eRe)$, then
$x,x'$ are isomorphic to $Re\otimes_{eRe}^{\mathbb{L}}y$,
$Re\otimes_{eRe}^{\mathbb{L}}y'$ for some $y,y'\in \D^c(eRe)$,
respectively. That is, the image of the functor
$Re\otimes_{eRe}^{\mathbb{L}}-: \D^c(eRe)\lra \D^c(R)$ is closed
under direct summands.

Indeed, let $y:=eR\otimes_R^{\mathbb{L}}x$ and
$y':=eR\otimes_R^{\mathbb{L}}x'$. Then it follows from
$$0=(R/J)\otimes_R^{\mathbb{L}}(Re\otimes_{eRe}^{\mathbb{L}}z)\simeq
(R/J)\otimes_R^{\mathbb{L}}(x\oplus x')\simeq
(R/J)\otimes_R^{\mathbb{L}}x\oplus (R/J)\otimes_R^{\mathbb{L}}x'$$
that $(R/J)\otimes_R^{\mathbb{L}}x\simeq 0\simeq
(R/J)\otimes_R^{\mathbb{L}}x'$. Now, by the definition of
recollements, there are two triangles in $D^c(R)$:
$$ Re\otimes_{eRe}^{\mathbb{L}}y\lra x\lra D(\lambda_*)\big((R/J)\otimes_R^{\mathbb{L}}x\big)\lra Re\otimes_{eRe}^{\mathbb{L}}y[1],$$
$$Re\otimes_{eRe}^{\mathbb{L}}y'\lra x'\lra
D(\lambda_*)\big((R/J)\otimes_R^{\mathbb{L}}x'\big)\lra
Re\otimes_{eRe}^{\mathbb{L}}y'[1].$$ Since the third terms of the
two triangles are isomorphic to zero, we get $x\simeq
Re\otimes_{eRe}^{\mathbb{L}}y \mbox{\; and \;} x'\simeq
Re\otimes_{eRe}^{\mathbb{L}}y'.$

By (1) and (2), we have verified all conditions of
Thomason-Waldhausen Localization Theorem \ref{waldhausen} for
$\D^c(R)$. This implies that all conditions of the
Thomason-Waldhausen Localization Theorem for $\Kb{\pmodcat{R}}$ are
satisfied, and therefore the sequence of the $K$-theory spaces:
$K(R/J)\longleftarrow K(R)\longleftarrow K(eRe)$ is a homotopy
fibration.


Let $\mathscr{P}^{<\infty}(R)$ be the full subcategory of $R$-mod
consisting of all $R$-modules $X$ with a finite-type resolution:
$$ 0\lra Q_m\lra \cdots\lra Q_1\lra Q_0\lra X\lra 0$$
for some $m\in \mathbb{N}$ such that all $Q_j$ are finitely
generated projective modules. By \cite[Section 4, Corollary
2]{Quillen}, we have $K(R)\simeq
K\big(\mathscr{P}^{<\infty}(R)\big)$. So, in the following we
identify $K(eRe)$ with $K\big(\mathscr{P}^{<\infty}(eRe)\big)$.

With this identification of $K$-theory spaces, now we show that the
map $K(Re\otimes_{eRe}-): K(eRe) \lra K(R)$ is a homotopy-split
injection, that is, there is a map $K(eR\otimes_R-): K(R)\ra K(eRe)$
of $K$-theory spaces, such that the composite of
$K(Re\otimes_{eRe}-)$ with $K(eR\otimes_R-)$  is homotopic to the
identity map on $K(eRe)$.

Consider the following commutative diagrams among exact categories
$\pmodcat{eRe}, \,\pmodcat{R}$ and $\mathscr{P}^{<\infty}(eRe):$
$$
\xymatrix{\,\pmodcat{eRe}\ar[r]^-{Re\otimes_{eRe}-}\ar@{_{(}->}[rd]&
\pmodcat{R}\ar[d]^-{eR\otimes_R-}\\
& \mathscr{P}^{<\infty}(eRe) } \Longrightarrow
\xymatrix{\,K(eRe)\ar[r]^-{K(Re\otimes_{eRe}-)}\ar[rd]^-{\simeq}_-{\mbox{homo.
equ}}&
K(R)\ar[d]^-{K(eR\otimes_R-)}\\
& K(\mathscr{P}^{<\infty}(eRe)) }\Longrightarrow
\xymatrix{\,K_n(eRe)\ar[r]^-{K_n(Re\otimes_{eRe}-)}\ar[rd]^-{\simeq}&
K_n(R)\ar[d]^-{K_n(eR\otimes_R-)}\\
& K_n(\mathscr{P}^{<\infty}(eRe)) }
$$
Note that the functor $eR\otimes_R-: \pmodcat{R}\lra
\mathscr{P}^{<\infty}(eRe) $ is well defined. Thus the map
$K(Re\otimes_{eRe}-): K(eRe)\lra K(R)$ is a homotopy-split
injection. By \cite{milnor2}, the $K$-theory spaces of rings are
always homotopy equivalent to CW-complexes. To conclude our
statement, we cite the following result in \cite[Corollary
7.15]{se}:

{\it For a homotopy fibration $X \lraf{f} Y \lra Z$ with $Z$
homotopy equivalent to a CW-complex, if  the map $f$ is
homotopy-split injection, then $Y$ is homotopy equivalent to the
product of $X$ and $Y$.}

\smallskip
Thus, from this result we know that the $K$-theory space $K(R)$ of
$R$ is homotopy equivalent to the product of the $K$-theory spaces
of $eRe$ and $R/J$, and therefore $$K_n(R)\simeq K_n(R/ReR)\oplus
K_n(eRe)$$ for all $n\in \mathbb{N}$. This completes the proof of
Proposition \ref{idempotent}, and also the proof of Theorem
\ref{thm1}. $\square$

\medskip
{\bf Proof of Theorem \ref{mainthm}}.

(1) Now, we assume that $\lambda :Y \ra X$ is a covariant morphism
of objects in ${\mathcal C}$. In this case, we consider
$\Lambda:=\End_{\mathcal C}(X\oplus Y)$, and let $J$ be the ideal of
$\Lambda$ generated by the projection $e$ from $X\oplus Y$ onto $Y$.
Then $e\Lambda e \simeq \End_{\mathcal C}(Y)$ and $\Lambda/J$ is
isomorphic to the quotient ring of $\End_{\mathcal C}(X)$ modulo the
ideal generated by those endomorphisms of $X$ which factorize
through the object $Y$, that is, $\Lambda/J\simeq \End_{{\mathcal
C},Y}(X)$ by Lemma \ref{lem01} (2). Since $_{\Lambda}J$ is
projective and finitely generated by Lemma \ref{lem01} (1), we can
apply Theorem \ref{thm1} to $\Lambda$ and $J$. In this case, we see
that the $K$-theory space of $\End_{\mathcal C}(X\oplus Y)$ is
homotopy equivalent to the product of the $K$-theory spaces of
$\End_{{\mathcal C},Y}(X)$ and $\End_{\mathcal C}(Y)$, and get
$$ K_n(\Lambda)\simeq K_n\big(\End_{{\mathcal C},Y}(X)\big)\oplus K_n\big(\End_{\mathcal C}(Y)\big)$$
for all $n\in\mathbb{N}$.

(2) Similarly, we may use Lemma \ref{3.4"} and Theorem \ref{thm1} to
show (2). $\square$

\medskip
{\it Remarks.} (1) If $I:=ReR$ is a finitely generated projective
$R$-module for $e^2=e\in R$, then $\End_R(I)$ and $eRe$ are Morita
equivalent. In fact, it follows from the projectivity of $_RI$ that
$Re\otimes_{eRe}eR\simeq ReR$ and that $eR$ is projective and
finitely generated as an $eRe$-module by Lemma \ref{projective}.
Thus we have $\add(Re)=\add(_RReR)$, where $\add(Re)$ stands for the
full subcategory of $\Modcat{R}$ consisting of all direct summands
of direct sums of finitely many copies of $Re$. This means that
$\End_R(_RI)$ is Morita equivalent to $eRe$.

(2) Dually, we may define \emph{weak trace factor modules}. Let $X$
and $Y$ be an $R$-modules. The module $Y$ is called a week trace
factor module of $X$ if there is a surjective homomorphism $\pi:
X\ra Y$ of $R$-modules such that the induced map $\pi^*:
\Hom_R(Y,Y)\lra \Hom_R(X,Y)$ is an isomorphism of abelian groups.
Obviously, in this case, the map $\pi$ is a contravariant
homomorphism.

\medskip
The dual of Theorem \ref{mainthm} can be stated as follows.

\begin{Theo}  Let $\mathcal C$ be an additive category and $f: Y\ra X$ a morphism of objects in
$\mathcal C$.

$(1)$ If $f$ is contravariant, then the $K$-theory space of
$\End_{\mathcal C}(X\oplus Y)$ is homotopy equivalent to the product
of the $K$-theory spaces of $\End_{\mathcal C}(X)$ and
$\End_{{\mathcal C},X}(Y)$. In particular, $$K_*\big(\End_{\mathcal
C}(X\oplus Y)\big)\simeq K_*\big(\End_{\mathcal C}(X)\big)\oplus
K_*\big(\End_{{\mathcal C}, X}(Y)\big)$$ for all $*\in \mathbb{N}$.

$(2)$ If $f$ is $Y$-contravariant, then the $K$-theory space of
$\End_{\mathcal C}(X\oplus Y)$ is homotopy equivalent to the product
of the $K$-theory spaces of $\End_{{\mathcal C}, Y}(X)$ and
$\End_{\mathcal C}(Y)$. In particular, $$K_*\big(\End_{\mathcal
C}(X\oplus Y)\big)\simeq K_*\big(\End_{{\mathcal C},Y}(X)\big)\oplus
K_*\big(\End_{\mathcal C}(Y)\big)$$ for all $*\in
\mathbb{N}$.\label{dualmainthm}
\end{Theo}

\medskip
The following is a consequence of Theorem \ref{mainthm} for
${\mathcal C}=R$-Mod.

\begin{Koro} Let $R$ be a ring with identity, and let
$\soc(R)$ be the socle of $_RR$. Then $$K_n\big(\End_R(R\oplus
\soc(R))\big)\simeq K_n\big(\End_R(\soc(R))\big)\oplus
K_n\big(R/\soc(R)\big)$$ for all $n\in \mathbb{N}.$ \label{socle}
\end{Koro}

{\it Proof.} Recall that for an $R$-module $M$, the socle of $M$ is
the sum of all simple submodules of $M$. Thus $\soc(R)$ is a direct
sum of minimal left ideals of $R$, and therefore it is actually an
ideal in $R$. Since $\soc(R)$ is a weak trace submodule of $_RR$ by
the definition of socles, we can apply Theorem \ref{mainthm} and get
$$K_n\big(\End_R(R\oplus \soc(R))\big)\simeq
K_n\big(\End_R(\soc(R))\big)\oplus K_n(R/\soc(R))$$ for all $n\in
\mathbb{N}$. $\square$

\medskip
For Auslander-Reiten sequences, we have the following result.

\begin{Koro} Let $A$ be an Artin algebra, and let $0\lra Z\lraf{g} Y\lraf{f} X\lra 0$ be
an Auslander-Reiten sequence in $A\modcat$.
If $\Hom_A(Y,Z)=0$, then

$$\begin{array}{rl}
K_n\big(\End_A(X\oplus Y)\big) & \simeq K_n\big(\End_A(Y)\big)\oplus K_n\big(\End_A(X)/\rad(\End_A(X))\big)  \\ &  \simeq
K_n\big(\End_A(Y)\big)\oplus K_n\big(\End_{A}(Z)/\rad(\End_A(Z))\big) \\ & \simeq K_n\big(\End_A(Y\oplus Z)\big)
\end{array} $$
for all $n\in
\mathbb{N}.$
\label{arsequence}
\end{Koro}

{\it Proof.} For Auslander-Reiten sequences, we know from \cite{hx2}
that $\End_A(Y\oplus Z)$ and $\End_A(X\oplus Y)$ are derived
equivalent, and therefore they have the isomorphic algebraic
$K$-groups, that is, $K_n(\End_A(X\oplus Y))\simeq K_n(\End_A(Y\oplus Z))$ for all $n\in \mathbb{N}.$
Note that $\Hom_A(Y,Z)=0$ if and only if the induced surjective map $\Hom_A(Y,Y)\ra \Hom_A(Y,X)$ is
an isomorphism of $\End_A(Y)$-modules. Since $f$ is surjective, it follows also from  $\Hom_A(Y,Z)=0$ that $\Hom_A(X,Z)=0$. Thus
$f: Y\ra X$ is a covariant homomorphism and, by
Theorem \ref{mainthm}, we have
$$K_n\big(\End_A(X\oplus Y)\big)\simeq K_n\big(\End_A(Y)\big)\oplus
K_n\big(\End_{A,Y}(X)\big)$$ for all $n\in \mathbb{N}.$

By properties of Auslander-Reiten sequences, we see that
$\rad(\End_A(X))$ is the image of the map $\Hom_A(X,f):
\Hom_A(X,Y)\ra \Hom_A(X,X)$. Thus $\End_{A,Y}(X)\simeq
\End_A(X)/\rad(\End_A(X))$ which is a division ring and isomorphic
to $\End_A(Z)/\rad(\End_A(Z))$. Hence
$$K_n\big(\End_A(X\oplus Y)\big)\simeq K_n\big(\End_A(Y\oplus Z)\big)\simeq K_n\big(\End_A(Y)\big)\oplus
K_n\big(\End_A(Z)/\rad(\End(Z))\big)$$ for all $n\in \mathbb{N}.$
$\square$

Further applications of Theorems \ref{thm1} and \ref{mainthm} will
be discussed in the next section.

\section{Applications\label{sect4}}

In this section, we deduce some consequences of Theorems \ref{thm1}
and \ref{mainthm}.

\subsection{Standardly stratified rings}

First, we consider standardly stratified rings and finite
dimensional quasi-hereditary algebras.

Standardly stratified and quasi-hereditary algebras were well
defined for finite dimensional algebras or semiprimary rings in
\cite{cps} and \cite{dr}, respectively. Now let us formulate them
for arbitrary rings.

Let $R$ be a ring with identity. Recall that an ideal $J$ of $R$ is
called a \emph{stratifying ideal }if

(1) $J=ReR$ for some idempotent element $e\in R$,

(2) $Re\otimes_{eRe}eR\simeq ReR$, and

(3) $\Tor_j^{eRe}(Re,eR)=0$ for $j\ge 1$.

\smallskip
Note that $J=ReR$ for $e^2=e\in R$ is a stratifying ideal if and
only if $J$ is homological. In particular, if $_RJ=ReR$ is
projective and finitely generated, then $J$ is a stratifying ideal
of $R$. In this case, $J$ is called a \emph{standardly stratifying}
ideal of $R$. The ring $R$ is called \emph{standardly stratified} if
there is a chain of ideals of $R$:
$$0=J_{n+1}\subseteq J_n\subseteq \cdots \subseteq J_2\subseteq J_1=R$$
such that $J_{i}/J_{i+1}$ is a standardly stratifying ideal in
$R/J_{i+1}$ for all $i$.

By this definition, every ring with identity is standardly
stratified. But the most interesting case for us is that for rings
we do have such a chain of length bigger than one.

For a standardly stratified ring $R$ with a defining chain of ideals
as above, there is an idempotent element $e\in R/J_{i+1}$ such that
$(R/J_{i+1})e(R/J_{i+1})= J_{i}/J_{i+1}$, we denote by $\Delta(i)$
the  $R$-module $(R/J_{i+1})e$. All these modules $\Delta(i)$ are
called the \emph{standard modules} with respect to the chain of
ideals of $R$. Note that standardly stratified rings may have
infinite global dimension.

By definition, a \emph{quasi-hereditary ring} is a standardly
stratified ring $R$ such that the endomorphism ring
$\End_R(\Delta(i))$ of each standard module $\Delta(i)$ is a
division ring. As in the case of finite dimensional algebras, one
can show that every quasi-hereditary ring has finite global
dimension.

Note that, by definition, the hereditary ring $\mathbb{Z}$ of
integers is a standardly stratified ring, but it is not a
quasi-hereditary ring. Thus, left hereditary rings may not be
quasi-hereditary. This example shows the difference of
quasi-hereditary rings defined in this paper from quasi-hereditary
algebras (or rings) in the sense of Cline, Parshall and Scott
\cite{cps} (or of Dlab and Ringel \cite{dr}). For further
information on finite dimensional standardly stratified and
quasi-hereditary algebras, we refer the reader to \cite{cps, dr} and
the references therein.

\begin{Koro} $(1)$ If $R$ is a standardly stratified
ring with the standard modules $\Delta(i)$ for $1\le i\le n$, then
$$K_*(R)\simeq \bigoplus_{j=1}^n K_*\big(\End_R(\Delta(j))\big)$$ for
all $*\in \mathbb{Z}$.

$(2)$ If $A$ is a finite-dimensional quasi-hereditary algebra over
an algebraically closed field $k$ with $n$ non-isomorphic simple
$A$-modules, then $G_*(A)\simeq K_*(A)\simeq nK_*(k)$ for all $*\in
\mathbb{N}$. \label{cor1}
\end{Koro}

{\it Proof.} (1) follows from Theorem \ref{thm1} inductively. (2) is
a consequence of (1) since for a finite-dimensional quasi-hereditary
algebra over an algebraically closed field $k$, we can refine a
hereditary chain into a maximal chain, and in this case, the
endomorphism ring of each standard module is isomorphic to the
ground field $k$.

Note that (2) follows also from the fact that finite-dimensional
quasi-hereditary algebras over an algebraically closed field $k$ are
noetherian and of finite global dimension. This implies that their
$G$-theory and $K$-theory  coincide. $\square$

\subsection{Matrix subrings}

In the following, we consider algebraic $K$-theory of matrix
subrings some of which are used in noncommutative algebraic geometry
(see \cite{ch2}) and arithmetic representation theory (see
\cite[Chapter 39]{cr}). In our discussions below, the key idea is to
find standardly stratifying ideals in those rings.

\begin{Koro} Let $R$ and $S$ be rings with identity, and let $_RM_S$ and $_SN_R$ be bimodule.
Suppose that $\varphi: M\otimes_SN\ra R$ and $\psi: N\otimes_RM\ra
S$ define a Morita context ring $T:=\begin{pmatrix}
  R             &  M    \\
 N    & S         \\
\end{pmatrix}$. If $\varphi$ is injective, and $_SN$ is projective and finitely
generated, then
$$K_n(T)\simeq K_n(S)\oplus K_n\big(R/(M\cdot N)\big)$$ for all $n\in \mathbb{N}$,
where $M\cdot N$ stands for the image of $\varphi$ in $R$,
\label{cor2'}
\end{Koro}

{\it Proof.} Let $e=\begin{pmatrix}
  0   &  0    \\
 0    & 1         \\
\end{pmatrix}$. Since $\varphi$ is injective and $_SN$ is projective and finitely
generated, it follows that $TeT=\begin{pmatrix}
  M\cdot N   &  M    \\
 N    & S         \\
\end{pmatrix}\simeq Te\otimes_SN \oplus Te$, which is a finitely generated projective
$T$-module, and that $T/TeT\simeq R/(M\cdot N)$. Thus, the corollary
follows immediately from Theorem \ref{thm1}. $\square$

\smallskip
{\it Remark.} The statement in Corollary \ref{cor2'} appeared in
\cite[Theorem 1.2]{Davydov}. However, the proof in \cite{Davydov}
seems to be wrong because the functor H in the proof is not well
defined.

\medskip
As a further corollary of Theorem \ref{thm1}, we consider the
question mentioned in Introduction (see also \cite{x1}) and provide
some partial answers. First of all, we mention the following
consequence of Corollary \ref{cor2'}, namely a result of Berrick and
Keating in \cite{bk}.

\begin{Koro} {\rm \cite{bk}} Let $R_i$ be a ring with identity for $i=1,2$,
and let $M$ be an $R_1$-$R_2$-bimodule. Then, for the triangular
matrix ring
$$ S=\begin{pmatrix}
  R_1             &  M    \\
  0    & R_2      &     \\
  \end{pmatrix},$$
there is an isomorphism of $K$-groups: $ K_n(S)\simeq K_n(R_1)\oplus
K_n(R_2)$ for $n\in\mathbb N$. \label{bktheorem}
\end{Koro}
In the next result, we consider slightly general matrix subrings.
Here, under the assumption that $_RJ$ is an idempotent, projective
and finitely generated ideal of a ring $R$, we extend the result
\cite[Proposition 5.3]{x1} for $K_1$ to a result for higher
$K$-groups.

\begin{Koro} Let $R$ be a ring with
identity, and let $J$ and $I_{i j}$ with $1\le i<j\le n$ be
arbitrary ideals of $R$ such that $I_{i j+1}J\subseteq I_{i\, j}$,
$JI_{i\, j}\subseteq I_{i+1 \,j}$ and $I_{i\, j}I_{j k}\subseteq
I_{i k}$ for $j<k\le n$. Define a ring
$$S:= {\begin{pmatrix}
  R & I_{1 2} &   \cdots &   \cdots & I_{1 n}\\
  J & R   &   \ddots  & \ddots &   \vdots\\
  J^2     & J  & \ddots  & \ddots  &   \vdots  \\
  \vdots  & \ddots &\ddots & R &  I_{n-1\, n}\\
  J^{n-1} &  \cdots & J^2 & J  &R\\
\end{pmatrix}.}_{n\times n}$$ If $_RJ$ is projective and finitely generated, then
$$K_*(S)\simeq  K_*(R)\oplus \bigoplus_{j=1}^{n-1}K_*(R/I_{j\;j+1}J).$$ \label{cor3}
\end{Koro}

{\it Proof.} We use induction on $n$ to prove this corollary.

Now let $e_i$ be the idempotent element of $S$ with $1_R$ at the
$(i,i)$-entry and zero at all other entries, and $e:=e_2+\cdots
+e_n$. As $J$ is a projective $R$-module, we have
$I_{ij}\otimes_RJ\simeq I_{ij}J$. Thus
$$SeS = {\begin{pmatrix}
  I_{12}J      & I_{12} &\cdots &\cdots    & I_{1n}\\
  J            & R      &\ddots &\cdots    & I_{2n}\\
  J^2            & J      & \ddots    &  \ddots &\vdots  \\
  \vdots       & \vdots &\ddots     & R & I_{n-1 n}\\
  J ^{n-1}           & J^{n-2}      &\cdots  & J   & R\\
\end{pmatrix}}\simeq Se\oplus Se_{2}\otimes_{e_2Se_2}J .$$
Here we identify $R$ with $e_2Se_2$. Since $_RJ$ is projective and
finitely generated, we infer that the $S$-module $SeS$ is also
projective and finitely generated. Clearly, $S/SeS$ is isomorphic to
$R/I_{12}J$. By Theorem \ref{thm1} (2), we get $K_*(S)\simeq
K_*(R/I_{12}J)\oplus K_*(eSe)$. By induction, we know that
$K_*(eSe)\simeq K_*(R)\oplus
\bigoplus_{j=2}^{n-1}K_*\big(R/I_{j\,j+1}J\big)$. Thus
$$K_*(S)\simeq  K_*(R)\oplus \bigoplus_{j=1}^{n-1}K_*(R/I_{j\;j+1}J).$$
This finishes the proof. $\square$

As a consequence of Corollary \ref{cor3}, we can prove the following
corollary.

\begin{Koro} Let $R$ be a ring with identity, and let $r$ be a
regular element of $R$ with $Rr=rR$. If $I_{i j}$ is an ideal of $R$
for $1\le i<j\le n$ such that $I_{i j+1}r\subseteq I_{i\, j}$,
$rI_{ij}\subseteq I_{i+1\,j}$ and $I_{i j}I_{j k}\subseteq I_{i k}$
for $j<k\le n$, then, for the matrix ring
$$T:=\begin{pmatrix}
  R & I_{12} &I_{13} & \cdots & I_{1\;n}\\
 Rr & R & I_{23} & \cdots & I_{2\;n} \\
\vdots & \ddots &\ddots &\ddots  & \vdots\\
 Rr^{n-2} & \cdots & Rr  & R& I_{n-1\;n}\\
Rr^{n-1}& \cdots &Rr^2 & Rr & R
  \end{pmatrix},$$ we have
$$K_*(T)\simeq K_*(R)\oplus \bigoplus_{j=1}^{n-1}K_*(R/I_{j\,j+1}r)$$ for all $n\in \mathbb{N}.$
\label{cor4}
\end{Koro}

By a \emph{regular element} we mean an element of $R$, which is not
a zero-divisor of $R$.

As an immediate consequence of Corollary \ref{cor4}, we mention the
following corollary for integral domain.

\begin{Koro} Let $D$ be an integral domain, $x\in D$,
and $I$ an ideal of $D$. Then
$$S:=\begin{pmatrix}
 D & I & \cdots & I\\
 Dx & \ddots & \ddots &\vdots \\
\vdots & \ddots &D   & I\\
 Dx^{n-1}& \cdots & Dx & D
  \end{pmatrix}_{n\times n}$$ is a ring, and
$$K_*(S)\simeq K_*(D)
 \oplus (n-1)K_*(D/Ix)$$ for all $*\in \mathbb{N}.$
 \label{cor5}
\end{Koro}

Now, we point out the following result.

\begin{Prop} Let $R$ be a commutative ring with identity, and let $x,y\in R$
such that $Rx+Ry=R$ and $Rx\cap Ry=Rxy$ (for example, $R$ is a
principle integral domain with $x$ and $y$ coprime in $R$). Suppose
that $y$ is invertible in an extension ring $R'$ of $R$. Then, for
the ring

$${S:=\begin{pmatrix}
 R & Rx & \cdots & Rx\\
 Ry & \ddots & \ddots &\vdots \\
\vdots & \ddots &R   & Rx\\
 Ry& \cdots & Ry & R
  \end{pmatrix},}_{n\times n}$$
we have $$ K_n(S)\simeq K_n(R)\oplus (n-1)K_n(R/Rx)\oplus
(n-1)K_n(R/Ry)$$ for all $n\in \mathbb N$.\label{4.7}
\end{Prop}

{\it Proof.} Let $\sigma$ be the diagonal matrix with the
$(1,1)$-entry $y$ and all other diagonal entries $1$. Then $\sigma$
is invertible in $M_n(R')$, the $n$ by $n$ full matrix ring of $R'$.
Let $B:= \sigma S\sigma^{-1}$. Thus $S\simeq B$ and $B$ is of the
form

$${B:=\begin{pmatrix}
 R     & Rxy   & Rxy   & \cdots & Rxy\\
 R     &  R    & Rx    & \cdots & Rx \\
 R     &  Ry   & R     & \ddots & \vdots\\
\vdots &\vdots &\ddots & R      &  Rx\\
 R     & Ry    &\cdots & Ry     &  R
  \end{pmatrix}.}_{n\times n}$$
Define $A:=M_n(R)$. Then $B$ is a subring of $A$ with the same
identity. Moreover, ${}_BA$ is isomorphic to the direct sum of $n$
copies of $Be_1$ where $e_1$ is the diagonal matrix
diag$(1,0,\cdots,0)$ of $B$. Thus $_BA$ is a finitely generated
projective $B$-module. Hence, by \cite[Lemma 3.1]{x1}, $B$ is
derived equivalent to $\End_B(B\oplus A/B)$. Clearly, the latter is
Morita equivalent to $\End_B(Be_1\oplus Q_2\oplus\cdots \oplus
Q_n)$, where $Q_j$ is given by the exact sequence
$$ 0\lra Be_j\lra Be_1\lra Q_j\lra 0, \;\;  2\le j\le n.$$
As in \cite[Section 3]{x1}, we can show that $\End_B(Be_1\oplus
Q_2\oplus\cdots \oplus Q_n)$ is isomorphic to the following ring
$$C:=\begin{pmatrix}
    R     & R/Rxy   & R/Rxy   & \cdots & R/Rxy\\
 0        &  R/Rxy  & Rx/Rxy  & \cdots & Rx/Rxy \\
 0        &  Ry/Rxy & R/Rxy   & \ddots & \vdots\\
\vdots    &\vdots   &\ddots   & \ddots  &  Rx/Rxy\\
 0        & Ry/Rxy  &\cdots   & Ry/Rxy &  R/Rxy
\end{pmatrix}.$$
From the Chinese remainder theorem we know  that $R/Rxy\simeq
R/Rx\oplus R/Ry$ as rings. Moreover, it follows from the assumptions
that the $R/Rxy$-bimodules $Rx/Rxy$ and $Ry/Rxy$ are isomorphic to
$R/Ry$ and $R/Rx$, respectively. Let $D$ be the lower right corner
$(n-1)\times (n-1)$-submatrix of $C$. Then the ring $D$ is actually
a direct sum of the following two rings:
$$ D={\begin{pmatrix}
 R/Ry   & R/Ry  & \cdots & R/Ry \\
 0  & R/Ry   & \ddots & \vdots\\
 \vdots  &\ddots   & R/Ry  &  R/Ry\\
 0  &\cdots   & 0 &  R/Ry
\end{pmatrix}_{n-1}}\bigoplus \quad {\begin{pmatrix}
 R/Rx   &  0 & \cdots & 0 \\
 R/Rx  & R/Rx   & \ddots & \vdots\\
 \vdots  &\ddots   & R/Ry  &  0\\
 R/Rx  &\cdots   & R/Rx &  R/Rx
\end{pmatrix}.}_{n-1}$$
Since derived equivalences preserve  algebraic $K_n$-groups (see
\cite{DS}), we have
$$ K_n(S)\simeq K_n(C)\simeq K_n(R)\oplus K_n(D)\simeq K_n(R)\oplus (n-1)K_n(R/Rx)\oplus (n-1)K_n(R/Ry) $$
for all $n\in \mathbb{N}.$ $\square$

\smallskip
{\it Remark.} For $n=2$, we can remove the conditions ``$Rx+Ry=R$
and $Rx\cap Ry=Rxy$" in Proposition \ref{4.7}, and get $K_*(S)\simeq
K_*(R)\oplus K_*(R/Rxy)$ for all $*\in \mathbb{N}$.

\medskip
Related to calculation of algebraic $K$-groups of the rings in the
proof of Proposition \ref{4.7}, the following result may be of
interest.

\begin{Koro} Let $R$ be a ring with identity, and let $I$ and $J$ be
ideals in $R$ with $JI=0$. If $_RI$ (or $J_R$) is projective and
finitely generated, then, for the ring

$${S:=\begin{pmatrix}
 R & I & \cdots & I\\
 J & R & \ddots &\vdots \\
\vdots & \ddots &\ddots   & I\\
 J& \cdots & J & R
  \end{pmatrix},}_{n\times n}$$
we have $$ K_*(S)\simeq nK_*(R)$$ for all $*\in \mathbb
N$.\label{4.8}
\end{Koro}

{\it Proof.} We assume that the $R$-module $_RI$ is projective and
finitely generated. Let $e:=e_1\in S$. Then
$$SeS:=
\begin{pmatrix}
 R & I & \cdots & I\\
 J & 0 & \cdots &0 \\
\vdots & \vdots & \ddots   & \vdots\\
 J& 0 & \cdots & 0
  \end{pmatrix}.$$
Since $_RI$ is projective, we have $J\otimes_RI \simeq JI = 0$ and
$Se\otimes_RI\simeq SeSe_j$ for $2\le j\le n$. Here we identify
$eSe$ with $R$. Since $_RI$ is projective and finitely generated, we
know that $_SSeSe_j$ is projective and finitely generated for $j=2,
\cdots, n$, and therefore the $S$-module $_SSeS\simeq Se\oplus
SeSe_2\oplus \cdots \oplus SeSe_n$ is a finitely generated
projective module. Thus, by Theorem \ref{thm1} and induction on $n$,
we have
$$ K_*(S)\simeq nK_*(R)$$
for all $*\in \mathbb{N}$.

The proof for the case that $J_R$ is projective and finitely
generated can be done similarly. $\square$

\medskip
{\it Remark.} If $R$ is an arbitrary ring with $I,J$ ideals in $R$
such that $IJ=JI=0$, then the ring $S$ in Corollary \ref{4.8} is the
trivial extension of $R\times R\times \cdots \times R$ by the
bimodule $L$, where $$L:=\begin{pmatrix}
 0 & I & \cdots & I\\
 J & 0 & \ddots & \vdots \\
\vdots & \ddots & \ddots   & I\\
 J& \cdots & J & 0
  \end{pmatrix}.$$
Thus we always have $K_n(S)\simeq nK_n(R)\oplus K_n(S,L)$ for all
$n\in \mathbb{N}$, where $K_n(S,L)$ is the $n$-th relative $K$-group
of $S$ with respect to the ideal $L$. This is due to the split
epimorphism $K_n(S)\ra K_n(S/L)$ of abelian groups, which is induced
from the split surjection $S\ra S/L$.

\medskip
Observe that rings of the form in Corollaries \ref{cor4}, \ref{cor5}
or Proposition \ref{4.7} occur in terminal orders over smooth
projective surfaces (see \cite{ch2}). For example, if we take $D$ to
be the power series ring $k[[z]]$ over a field $k$ in one variable
$z$, $I=zk[[z]]$ and $x=1$, then the ring $S$ in Corollary
\ref{cor5} is related to the completion of a closed point in a
quasi-projective surface.  It would be interesting to know how
$K$-theory or recollements could be applied in this situation.

\subsection{Some special rings \label{sect4.2}}

In this section, we consider the algebraic $K$-theory of rings
appearing in different areas.

\subsubsection{Algebraic $K$-theory for affine cellular algebras}

As a generalization of cellular algebras in the sense of Graham and
Lehrer \cite{gl}, affine cellular algebras were introduced in
\cite{kx} to study the representation theory and homological
properties of certain infinite dimensional algebras which include
extended affine Hecke algebras of type $\tilde{A}$. We shall see
that the $K$-theory of affine cellular algebras can be studied in
local information. First, we recall the definition of affine
cellular algebras from \cite{kx}.

Let $k$ be a noetherian integral domain. For two $k$-modules $W$ and
$V$, we denote the switch map by $\tau:   W\otimes_kV\ra
V\otimes_KW, w\otimes v\mapsto v\otimes w$ with $w\in W$ and $v\in
V$.

\begin{Def}  {\rm \cite{kx}} Let $A$ be a unitary $k$-algebra with a
$k$-involution $i$ on $A$. A two-sided ideal $J$ in $A$ is called an
{\em affine cell ideal} if and only if the following data are given
and the following three conditions are satisfied:
\begin{enumerate}
\item[$(1)$]
The ideal $J$ is fixed by $i$, that is, $i(J) = J$.

\item[$(2)$] There exist a free $k$-module $V$ of
finite rank and a finitely generated commutative  $k$-algebra $B$
with identity and with a $k$-involution $\sigma$ such that $\Delta
:=V\otimes_kB$ is an $A$-$B$-bimodule, where the right $B$-module
structure is induced by that of the right regular $B$-module $B_B$.

\item[$(3)$] There is an $A$-$A$-bimodule isomorphism $\alpha:
J\lra \Delta\otimes_B\Delta'$, where $\Delta' =B\otimes_kV$ is a
$B$-$A$-bimodule with the left $B$-structure induced by $_BB$ and
with the right $A$-structure via $i$, that is, $(b\otimes v)a :=
\tau(i(a)(v\otimes b))$ for $a\in A, b\in B$ and $v\in V$), such
that the following diagram is commutative:

$$\begin{array}{rll}
 J &\stackrel{\alpha}{\longrightarrow}&\Delta\otimes_B \Delta'\\
i\Big\downarrow &         &\quad\Big\downarrow{v_1\otimes
b_1\otimes_B
b_2\otimes v_2\mapsto v_2\otimes \sigma(b_2)\otimes_B \sigma(b_1)\otimes v_1}\\
 J &\stackrel{\alpha}{\longrightarrow}&\Delta\otimes_B \Delta'
\end{array}$$
\end{enumerate}

The algebra $A$ (with the involution $i$) is called an {\em affine
cellular algebra} if and only if there is a $k$-module decomposition
$A=J_1' \oplus J_2' \oplus \dots \oplus J_n'$ (for some $n$) with
$i(J_j')=J_j'$ for each $j$, such that setting $J_j = \oplus_{l=1}^j
J_l'$ gives a chain of two-sided ideals of $A$:
$$0=J_0 \subset J_1 \subset J_2 \subset \dots \subset J_n = A$$ (each
of them fixed by $i$) and for each $j$ ($j = 1, \dots, n$) the
quotient $J_j'=J_j/J_{j-1}$ is an affine cell ideal of $A/J_{j-1}$
(with respect to the involution induced by $i$ on the quotient).
\label{ourdef}
\end{Def}
By definition, for each subquotient $J_j/J_{j-1}$ of an affine
cellular algebra $A$, there is a commutative algebra $B_j$  and an
$A$-$B_j$-bimodule $\Delta(j)$ such that $J_j/J_{j-1}$ is an affine
cell ideal in $A/J_{j-1}$. In this case, we say that $B_j$ is
\emph{associated} with $J_j/J_{j-1}$, and $\Delta(j)$ is a
\emph{cell module}.

\begin{Prop} Let $A$ be an affine cellular algebra with a cell chain
$J_0=0\subset J_1\subset \dots \subset J_n=A$ and the associated
commutative rings $B_j$ for $1\le j\le n$. Suppose that each $B_j$
satisfies rad$(B_j)=0$ and that each $J_j/J_{j-1}$ is idempotent and
contains a non-zero idempotent element in $A/J_{j-1}$. Then
$$K_*(A)\simeq \displaystyle\bigoplus_{j=1}^nK_*(B_j)$$ for all $*\in \mathbb{N}$.
\label{affcelalg}
\end{Prop}

{\it Proof.} Under the assumptions of Proposition \ref{affcelalg},
we know from the proof of \cite[Theorem 4.3]{kx} that each ideal
$J_j/J_{j-1}$ of $A/J_{j-1}$ is generated by an idempotent element
$e_j$ and that, as an $A/J_{j-1}$-module, $J_j/J_{j-1}$ is
projective and isomorphic to a direct sum of finitely many copies of
the cell module $\Delta(j)$. Moreover, it follows from the proof of
\cite[Theorem 4.3]{kx} that
$\add(_A\Delta(j))=\add(_A(A/J_{j-1})e)$. This implies that
$\Delta(j)$ is a finitely generated $A$-module and that
$e_j(A/J_{j-1})e_j$ is Morita equivalent to $B_j$. Thus we may
inductively apply Theorem \ref{thm1} (2) to get Proposition
\ref{affcelalg}. $\square$

\subsubsection{Algebraic $K$-theory for affine Hecke algebras and
quantum Schur algebras}

Let $k$ be the Laurent polynomial ring $\mathbb{Z}[q,q^{-1}]$ in
variable $q$ over the ring $\mathbb{Z}$ of integers. Let $(W, S)$ be
a Coxeter system. For example, if $W$ is the symmetric group on the
letters $\{1, 2, \cdots, n\}$ with $S:=\{s_i=(i,i+1)\mid i=1,2,
\cdots, n-1\}\subseteq W$, then the Coxeter system is said to be of
type $A_{n-1}$. The \emph{Hecke algebra} of $(W,S)$ over $k$,
denoted by $\mathcal{H}_k(W,S)$, is a unitary associative algebra
with a $k$-basis $\{T_w\mid w\in W\}$, subject to the following
relations:

$$\begin{array}{ll} (T_s-q^2)(T_s+1)=0\quad & \mbox{if}\; s\in S,\\
T_wT_u=T_{wu} \quad & \mbox{if}\; \ell(wu) = \ell(w)+\ell(u),
\end{array}$$
where $\ell$ is the usual length function of $W$.

Let $(W,S)$ be the Coxeter system of type $\tilde{A}_{n-1}$. Then
the cyclic group $\mathbb{Z}/n\mathbb{Z}$ acts on $W$. Thus we may
form the semiproduct $\widetilde{W}:=W\ltimes {\mathbb{Z}}/n{\mathbb
Z}$, and define similarly the Hecke algebra over $k$ of the extended
Coxeter system $(\widetilde{W}, S)$. This Hecke algebra is then
called the \emph{extended affine Hecke algebra} of type
$\tilde{A}_{n-1}$, denoted by $\h_k(n,r)$. For more details about
affine Hecke algebras we refer to \cite{lusztig}.

We may apply Proposition \ref{affcelalg} to the extended affine
Hecke algebra of type $\widetilde{A}_n$ since this algebra was shown
to be affine cellular in \cite{kx}. The proofs there imply the
following corollary.

\begin{Koro} Let $k$ be a field of characteristic $0$ and $q\in k$ such that $\sum_{w\in W_0}q^{\ell(w)}\ne 0$,
where $W_0$ is the symmetric group of $n$ letters. For the extended
affine Hecke algebra $\h_k(n,r)$, we have
$$ K_*(\h_k(n,r)) \simeq \bigoplus_{\bf c}K_*(R_{{\bf c}})$$
where $c$ runs over all two-sided cells of the extended affine Weyl
groups $\tilde{W}$, and $R_{\bf c}$ stands for the representation
ring associated with ${\bf c}$, which is isomorphic to a tensor
product of rings of the form $\mathbb{Z}[X_1,
X_2,\cdots,X_{s+1}]/(X_sX_{s+1}-1)$. \label{affine}
\end{Koro}

Now, we turn to quantum Schur algebras. Let $(W,S)$ be the Coxeter
system of type $A_{n-1}$, and $\mathcal{H}_q(n)$ be its Hecke
algebra over $k$. Given a partition $\lambda$ of $n$, one may define
a Young subgroup $W_{\lambda}$ of $W$, and an element
$x_{\lambda}:=\sum_{w\in W_{\lambda}}q^{\ell(w)}T_w$. Suppose $r\le
n$. Let $\Lambda^+(n,r)$ be the set of partitions of $n$ with at
most $r$ parts. The \emph{quantum Schur algebra}
$\mathcal{S}_q(n,r)$ is defined as
$$ \mathcal{S}_q(n,r):=\End_{\mathcal{H}_q(n)}\big(\bigoplus_{\lambda\in\Lambda^+(n,r)}\mathcal{H}_q(n)x_{\lambda}\big).$$
Quantum Schur algebras have many nice homological properties, for
example, they are (integral) quasi-hereditary algebras over $k$ and
their standard modules $\Delta(\lambda)$, indexed by
$\Lambda^+(n,r)$, have the property that
$\End_{S_q(n,r)}(\Delta(\lambda))\simeq k$ for all $\lambda\in
\Lambda^+(n,r)$. Thus, by Corollary \ref{cor1}, we have the
following result.

\begin{Koro} For the quantum Schur algebra $\mathcal{S}_q(n,r)$, we have
$$ K_*(\mathcal{S}_q(n,r)) \simeq mK_*(\mathbb{Z})\oplus mK_{*-1}(\mathbb{Z})$$
for all $*\in \mathbb{N}$, where $m$ is the cardinality of the set
$\Lambda^+(n,r).$ \label{qschur}
\end{Koro}

{\it Proof.} Since the ring $\mathbb{Z}$ is noetherian and of finite
global dimension, we know that $$K_i(\mathbb{Z}[t, t^{-1}])\simeq
K_i(\mathbb{Z})\oplus K_{i-1}(\mathbb{Z})$$(see \cite[Theorem
8]{Quillen}). 
Thus Corollary \ref{qschur} follows immediately from Corollary
\ref{cor1}. $\square$

\subsection{Algebraic $K$-theory for skew group rings\label{sect4.3}}

Let $S$ be a ring with identity and suppose that $G$ is a finite
group of automorphisms of the ring $S$ such that the order of $G$ is
a unit in $S$. Let $R= S*G$ be the skew group ring and
$e:=\frac{1}{|G|}\sum_{g\in G}g$. Then $e^2=e$ and the ring
$S^G:=\{s\in S\mid s^g=s \mbox{ for all }g\in G\}$ of invariants of
$G$ is isomorphic to $\End_R(Re)$. We write $R(S,G)$ for the
ramification algebra $R/ReR$. The trace ideal of $Re$ in $R$ is
$ReR$. If $_RReR$ (or $ReR_R$) is projective and finitely generated,
then $\End_R(ReR)$ is Morita equivalent to $S^G$. Thus we have the
following consequence of Theorem \ref{thm1}.

\begin{Koro} For any $n\in \mathbb{N}$, there hold:

$(1)$ $K_n\big(\End_R(R\oplus ReR)\big)\simeq
K_n\big(R(S,G)\big)\oplus K_n\big(\End_R(ReR)\big)$.

$(2)$ If the $R$-module $_RReR$ or $ReR_R$ is finitely generated and
projective, then $$K_n(R)\simeq K_n\big(R(S,G)\big) \oplus
K_n(S^G).$$ \label{invariants}
\end{Koro}

Note that the case $ReR=R$ was considered in \cite{lorenz} to
compare the $K_0$-groups of $R$ with those of $S$ and $S^G$. In
fact, in this case, the condition $|G|^{-1}\in S$ implies that $R$
and $S^G$ are Morita equivalent. The higher algebraic $K$-groups of
$S^G$ were discussed in \cite{hodges} under some additional
assumptions on both $S$ and $G$.

Let us mention an example in which the second condition in Corollary
\ref{invariants} holds true. For instance, if $S$ is a finite
product of simple rings, and if $G$ is a finite group acting as
automorphisms of $S$ such that the order of $G$ is invertible in
$S$, then the skew group ring $R$ is also a finite product of simple
rings. Thus $_RReR$ is projective.

\section{Examples\label{example}}
The following examples illustrate how our results in this note can
be used to compute algebraic $K_n$-groups of rings. They also show
that some conditions on $I$ in Theorem \ref{thm1} cannot be omitted
or weakened.

\medskip
{\bf Example 1.} Let $k$ be a field, and let $R$ be the ring
$k[X]/(X^2)$. We denote by $x$ the element $X +(X^2)$ in $R$. Then
we may form the matrix ring
$$ A:={\begin{pmatrix}
R & k  \\
k& k \end{pmatrix}}, \quad {\begin{pmatrix}
r+sx & a  \\
b& c \end{pmatrix}\begin{pmatrix}
r'+s'x & a'  \\
b'& c' \end{pmatrix}=\begin{pmatrix}
rr'+(rs'+sr'+ab')x & ra'+ac'  \\
br'+cb'& cc' \end{pmatrix}}$$ for $r,r',s,s',a,a',b,b',c,c'\in k$.
One can check that this matrix ring is isomorphic to the quotient
algebra of the path algebra of the quiver
$$ 
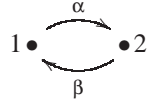
\xymatrix{1\,\bullet\ar^-{\bf{\alpha}}@/^0.8pc/[r]
&\bullet\, 2\ar^-{\bf{\beta}}@/^0.8pc/[l]}$$ modulo the ideal
generated by $\beta\alpha$, and the latter is
is a quasi-hereditary $k$-algebra of global dimension $2$. Note also
that $A\simeq \End_R(R\oplus k)$ and that $K_*(A)\simeq K_*(k)\oplus
K_*(k)$ by Theorem \ref{thm1}. If $k$ is a finite field, then we
have a full knowledge of $K_*(A)$ by a result in \cite{qff}.

Let $I=AeA$, where $e$ is the idempotent of $A$ corresponding to the
vertex $1$. Then $I$ is an idempotent ideal in $A$, and $_AI$ is
finitely generated, and has finite projective dimension, but not
projective. Clearly, we have $K_n(eAe)\simeq K_n(R)$ and
$K_n(A/AeA)\simeq K_n(k)$. Clearly, $K_0(A)\not\simeq K_0(A/I)\oplus
K_0\big(\End_A(I)\big)$. Hence, if $_AI$ is not projective, then the
second statement of Theorem \ref{thm1} (1) may fail in general.
Since $K_1(R)\simeq k\oplus \; k^{\times}$, we get $K_1(A)\not\simeq
K_1(eAe)\oplus K_1(A/AeA)$. Note that in this example, the condition
that $\Tor_j^R(A/I,A/I)=0$ for all $j>0$ fails, that is, the ideal
$I$ is not homological. Thus, in Theorem \ref{thm1} (2), that $_RI$
is homological cannot be dropped. This example also shows that
Corollary \ref{cor2'} may be false if $_SN$ is projective but
$\varphi$ is not injective.

If we modify this example slightly and just consider the algebra $B$
given by the above quiver but with the relation
$\beta\alpha\beta=0$, then the ideal $I'=Be_1B$ is homological with
infinite projective dimension as a left $B$-module. In this case,
$K_*(B)\simeq K_*(R)\oplus K_*(k)$ by Corollary \ref{cor1}(1). But,
since $\End_B(I')\simeq A$ and $B/I'\simeq k$, we cannot get
$K_n(B)\simeq K_n(B/I')\oplus K_n(\End_B(I'))$.  This shows that the
projectivity of $I$ in the second statement of Theorem \ref{thm1}
(1) cannot be relaxed to homological ideal.

{\bf Example 2.} Note that for the  triangular matrix ring
$T:=\begin{pmatrix}
k & k  \\
0& k \end{pmatrix}$, if we take $I:=\begin{pmatrix}
0 & k  \\
0& 0 \end{pmatrix}$, then $I^2=0$ and $_TI$ is projective and
finitely generated. In this case, we can see that $K_n(T)\simeq
K_n(k)\oplus K_n(k)$ by Corollary \ref{bktheorem}. Thus
$K_n(T/I)\oplus K_n(\End_T(I))=K_n(k)\oplus K_n(k)\oplus
K_n(k)\not\simeq K_n(T) $. Hence the condition $I^2=I$ in Theorem
\ref{thm1} (1) cannot be removed.

\medskip
{\bf Example 3.} Let $p>0$ be a prime integer, and let
$\mathbb{Z}_p$ be the ring of $p$-adic-integers.  We consider the
ring
$$ R:=\begin{pmatrix}
\mathbb{Z}_p & p\mathbb{Z}_p  \\
p\mathbb{Z}_p& \mathbb{Z}_p\end{pmatrix}. $$ Clearly,
$p\mathbb{Z}_p$ is isomorphic to $\mathbb{Z}_p$ as left
$\mathbb{Z}_p$-modules, and therefore is projective and finitely
generated. Thus, by Corollary \ref{cor5} (see also Remark to
Proposition \ref{4.7}), we have
$$ K_n(R)\simeq K_n(\mathbb{Z}_p )\oplus
K_n(\mathbb{Z}_p/p^2\mathbb{Z}_p)\simeq K_n(\mathbb{Z}_p )\oplus
K_n(\mathbb{Z}/p^2\mathbb{Z}).$$

\medskip
Based on the results and examples in this note, we mention the
following questions.

\medskip
{\bf Open questions.}  (1) Let $R$ be a ring with identity and $I$
an ideal of $R$ with $I^2=0$. We define a ring $S:=\begin{pmatrix}
R & I  \\
I & R\end{pmatrix}.$ How is the algebraic $K$-group $K_n(S)$ of $S$
related to the $K_n$-groups of rings produced from $R$ and the ideal
$I$ for $n\ge 2$?

\medskip
Note that $K_i(S)=K_i(R)\oplus K_i(R)$ for $i=0,1.$ This can be done
by using Mayer-Vietoris sequences.

\medskip
(2) Let $R$ be a ring with identity and $e=e^2\in R$. Suppose that
$ReR$ is homological and $_RReR$ possesses an infinite resolution by
finitely generated projective $R$-modules. Does the following
isomorphism hold true:
$$ K_n(R)\simeq K_n(R/ReR)\oplus K_n(eRe)$$
for every $n\in \mathbb{N}$?

\smallskip
In the question (2), the canonical surjection $R\ra R/ReR$ is a ring epimorphism.
From the representation-theoretic point of view, the ring $R/ReR$ has less simple modules than $R$ does. 
There are ring epimorphisms for which the cardinality of simple modules may increase (see \cite{xc2}). 
They arise from the so-called universal localizations. 
In another paper we shall establish a formula for higher algebraic 
$K$-groups of two rings linked by such a ring epimorphism.

\bigskip
{\bf Acknowledgements.} The corresponding author C.C. Xi is grateful
to Nanqing Ding, Xuejun Guo and Yuxian Geng for some discussions and
their hospitality during a visit to Nanjing University and Jiangsu
Teachers University of Technology in April, 2012. The research work
of C.C. Xi is partially supported by DPF(20110003110003) and PCSIRT,
while the one of H.X. Chen is partially supported by Post-Doctor
Funds of China.

Part of the paper is revised during Xi's visit to the Institute of
Algebra and Number Theory, University of Stuttgart, in July and
August, 2012. He would like to thank Steffen K\"onig for the
invitation and his hospitality.

\medskip

{\footnotesize

\begin{tabular}{lcl}
Hongxing Chen, & \qquad \qquad\qquad & Changchang Xi \\
Beijing International Center for Mathematical Research, & \qquad &
School of
Mathematical Sciences\\
Peking University & \qquad & Capital Normal University\\
100871 Beijing & \qquad & 100048 Beijing \\
People's Republic of China & \qquad & People's Republic of China\\
{\tt Email: chx19830818@163.com} & \qquad & {\tt Email:
xicc@bnu.edu.cn}
\end{tabular}

\bigskip
First version: November 16, 2011, revised: July 23, 2012}
\end{document}